\documentclass[11pt,reqno]{amsart}
\usepackage{pb-diagram} 
\usepackage{amsmath,amssymb}

\newtheorem{definition}{Definition}
\newtheorem{theorem}{Theorem}
\newtheorem{proposition}[theorem]{Proposition}
\newtheorem{corollary}[theorem]{Corollary}
\newtheorem{remark}{Remark}
\newtheorem{lemma}[theorem]{Lemma}
\newfont{\bb}{msbm10 at 12pt}

\def\eucl{\mathrm{eucl}}

\def\Si{{\Sigma}}

\def\Ga{{\Gamma}}

\def\phi{{\varphi}}

\def\TM{{\mathrm{T}\mathrm{M}}}  
   
\def\M{\mathrm{M}}
\def\T*M{\mathrm{T}^*\mathrm{M}}

\def\H{\mathrm{H}}
\def\R{\mathrm{R}}
\def\gst{{g_{\rm st}}}

\def\L{\mathrm{L}}

\def\MIT{\mathrm{MIT}}
\def\CHI{\mathrm{CHI}}
\def\L{\mathrm{L}}
\def\BMIT{\mathbb{B}_{\MIT}}
\def\BCHI{\mathbb{B}_{\CHI}}

\def\D{\mathrm{D}}
\def\hs{\mathbb{S}^n_+}

\def\V{\mathrm{V}}
\def\G{\mathrm{G}}
\def\ov{\overline}
\def\F{\mathrm{F}}
\def\U{\mathrm{U}}
\let\pa\partial     
\let\na\nabla     
     
\DeclareMathAlphabet{\doba}{U}{msb}{m}{n}

\def\T{\mathrm{T}}
\def\Vol{{\mathop{\rm Vol}}}

\def\Vol{{\mathop{\rm Vol}}}

\def\BMIT{\mathbb{B}_{\MIT}} 
\def\BCHIpm{\mathbb{B}_{\CHI}^\pm}

\def\bhlp{\lambda_{\min}^+(\M,[g],\sigma)}

\def\Ker{\mathrm{Ker}}
\def\Coker{\mathcal{C}}
\let\<\langle     
\let\>\rangle

\def\Id{{\mathop{\rm Id}}}

\def\hpm{\mathcal{H}^{\pm}}

\let\<\langle     
\let\>\rangle

\long\def\komment#1{} 

\begin{document}


\title[Green functions and applications]{Green functions for the Dirac operator under local boundary conditions and applications}     
\author{Simon Raulot} 
\address{Institut {\'E}lie Cartan\\
Universit{\'e} Henri Poincar{\'e}, Nancy I\\
B.P. 239\\
 54506 Vand\oe uvre-L{\`e}s-Nancy Cedex, France}
\email{raulot@iecn.u-nancy.fr}  
\date{\today}


\maketitle

\begin{abstract}
In this paper, we define the Green function for the Dirac operator under two local boundary conditions: the condition associated with a chirality operator (also called the chiral bag boundary condition) and the $\MIT$ bag boundary condition. Then we give some applications of these constructions for each Green function. From the existence of the chiral Green function, we derive an inequality on a spin conformal invariant which, in particular, solve the Yamabe problem on manifolds with boundary in some cases. Finally, using the $\MIT$ Green function, we give a simple proof of a positive mass theorem previously proved by Escobar. 
\end{abstract}


\section{Introduction}


A nice and powerful application in the study of the Green function of an elliptic differential operator appears in Schoen's approach \cite{schoen:84} of the Yamabe problem. This famous problem of conformal Riemannian geometry can be stated as follow: given a connected compact Riemannian manifold $(\M^n,g)$ with $n\geq 3$, can one find a metric $\ov{g}$ in the conformal class $[g]$ of $g$ with constant scalar curvature? We first briefly recall some steps in the resolution of this problem. In \cite{yamabe:60}, Yamabe claims to answer this question, however Tr\"udinger \cite{trudinger:68} points out a mistake in Yamabe's argument and can fix the proof in some cases. The remaining cases were treated by Aubin \cite{aubin:76} and Schoen \cite{schoen:84}. In fact, solving the Yamabe problem is equivalent to find a smooth positive solution of a non linear elliptic equation involving a conformally invariant operator: the conformal Laplacian. This is a two order elliptic differential operator acting on functions defined by:
\begin{eqnarray*}
\L_g :=\frac{4(n-1)}{n-2}\Delta_g +\R_g 
\end{eqnarray*}

and which relates the scalar curvature of two metrics in the same conformal class by the formula:
\begin{eqnarray*}
\R_{\ov{g}}=f^{-\frac{n+2}{n-2}}\L_g f,
\end{eqnarray*}

if $\ov{g}=f^{\frac{4}{n-2}}g\in[g]$ and where $\Delta_g$ and $\R_g$ denote respectively the standard Laplacian and the scalar curvature of $(\M,g)$. Then finding a metric with constant scalar curvature in the conformal class of $g$ is equivalent to find a smooth positive function $f$ for the equation:
\begin{eqnarray}\label{equationyam}
\L_g f= {\rm{C}}\, f^{N-1}
\end{eqnarray}

where $\rm{C}$ is a constant and $N=\frac{2n}{n-2}$. The major difficulty in this problem comes from the fact that the Sobolev embedding $\rm{H}^2_1(\M)\hookrightarrow\rm{L}^N(\M)$ is not compact and then a standard variationnal approach cannot allow to conclude. However with a little work, one can show the existence of a non negative smooth function $f$ solution of Equation~(\ref{equationyam}) which is either positive or zero on $\M$. So it remains to show that this function cannot vanish identically on $\M$. In \cite{aubin:76}, Aubin observed that if:
\begin{eqnarray}\label{aubin}
\mu(\M)<\mu(\mathbb{S}^n)
\end{eqnarray}
 
then the function $f$ has to be positive. Here $\mu(\M)$ denotes the Yamabe invariant of $\M$ (which only depends of the conformal class of $g$) defined by:
\begin{eqnarray}\label{yamcomp}
\mu(\M):=\underset{u\in\mathrm{C}^1(\M),u\neq 0}{\mathrm{inf}} \frac{\int_{\M}u{\rm{L}}_g u\, dv(g)}{\Big(\int_{\M}|u|^N dv(g)\Big)^{\frac{2}{N}}}.
\end{eqnarray}

He also proved that if $(\M,g)$ is not conformally equivalent to the standard sphere $(\mathbb{S}^n,g_{\rm{st}})$ then $\mu(\M)\leq\mu(\mathbb{S}^n)$ and moreover if $n\geq 6$ and the manifold is not locally conformally flat then (\ref{aubin}) holds. The last step in the resolution of the Yamabe problem was made by Schoen \cite{schoen:84} which proved $(\ref{aubin})$ in the remaining cases using the positive mass theorem. In $1987$, Lee and Parker \cite{lee.parker:87} gave a proof of the Yamabe problem unifying Aubin's and Schoen's arguments. The key point of their work uses Schoen's result and is strongly based on the study of the Green function of the conformal Laplacian and on the existence of an adaptated coordinates chart namely the conformal normal coordinates. In fact, using the development of the Green function near a point $q\in\M$ in these coordinates, they construct adapted test functions which include Aubin and Schoen's arguments. More precisely, in a neighbourhood of $q$, the conformal Laplacian's Green function can be decomposed into the sum of a singular part (which is nearly the Green function for the standard Laplacian in Euclidean space) and a regular part. The fact is that according to the manifold is or is not conformally flat, this regular part evaluated in $q$ is either the mass of an asymptotically flat manifold (obtained by conformal change of metrics whose weight is exactly the Green function) either the squared norm of the Weyl tensor at $q$.\\

If now we suppose that the manifold $\M$ is spin, there exists another operator which has the same conformal properties than the conformal Laplacian operator: the Dirac operator. This operator is a first order elliptic differential operator acting on sections of the complex spinor bundle over $(\M,g)$. The relation between these two operators was made by Hijazi (see \cite{hijazi:86}) with an inequality (now called the Hijazi inequality) which relates their first eigenvalues and which allows to define a spin conformal invariant given by: 
\begin{eqnarray*}
\bhlp:=\underset{\ov{g}\in[g]}{\inf}\big\{\lambda_1^+(\ov{g})\Vol(\M,\ov{g})^{\frac{1}{n}}\big\}.
\end{eqnarray*}

This invariant can be seen as an analogue of the Yamabe invariant in the spinorial setting. Moreover, using the Hijazi inequality, we can compare this invariant to the Yamabe invariant of $(\M,g)$ (see \cite{hijazi:91} and \cite{baer:92} for the case $n=2$):
\begin{eqnarray}\label{hijazicomp}
\bhlp^2 \geq\frac{n}{4(n-1)}\mu(\M).
\end{eqnarray}

This invariant has been (and is again) the main subject of many works, in particular in a serie of papers by Ammann, Humbert and Morel (see \cite{ammann}, \cite{amm1}, \cite{amm4} or \cite{ahm}). In these papers, the authors show that this invariant shares a lot of the Yamabe's invariant properties. Indeed, in \cite{amm3} it is shown that we can derive a spinorial analogue of Aubin's inequality, i.e. we have:
\begin{eqnarray}\label{ahm}
\bhlp\leq \lambda_{\min}(\mathbb{S}^n,[\gst],\sigma_{\mathrm{st}})=\frac{n}{2}\omega_n^{\frac{1}{n}}
\end{eqnarray}

where $\omega_n$ denotes the volume of the standard sphere. Moreover, if this inequality is strict, then with the help of the Hijazi inequality, it gives a spinorial proof of the Yamabe problem. On the other hand, Ammann \cite{habilbernd} also proves that the strict inequality in (\ref{ahm}) gives the existence of a nonlinear partial differential equations for the Dirac operator involving a critical Sobolev embedding theorem. A natural question is then to ask if we can find sufficient conditions for which Inequality~(\ref{ahm}) is strict. A partial answer is given in \cite{amm4} with the help of a detailed study of the Green function for the Dirac operator. Indeed, in this article, they considered the case of locally conformally flat manifolds where the behaviour of the Green function can be well understood. In fact, they show that the Dirac operator's Green function can be decomposed (around $q$) into the sum of a singular part (which is nearly the Green function for the Dirac operator in the Euclidean space) and a smooth spinor field. Then they define the mass operator which plays the same role as the constant term of the Green's function of the conformal Laplacian. An other application of the Green function of the Dirac operator can also be found in \cite{amm2} where a simple proof of the positive mass theorem for the Yamabe problem is given. The main idea of the proof is based on the fact that the Witten's spinor (see \cite{lee.parker:87} or \cite{bartnik}) is exactly obtained as the image of the Green function's of the Dirac operator under the conformal change of metric whose weight is given by the Green function of the conformal Laplacian.\\

In this paper, we define and study the Green's function of the Dirac operator on manifolds with boundary. More precisely, since the conformal aspect is underlying in this work, the choice of the boundary conditions is crucial for the following. It turns out that local boundary conditions seems to be appropriate for the study of the conformal aspect of the Dirac operator on manifolds with boundary.\\

In a first part, we will focus on the condition associated with a chirality operator (also called the chiral bag boundary condition, see \cite{hijazi.montiel.roldan:01} for example). A direct application of this construction is motivated by previous results of the author (see \cite{sr}, \cite{sr3} or Section~\ref{apl1}). Indeed, in \cite{sr3} we define an analogue of $\bhlp$ on manifolds with boundary and we show that this invariant satisfies an inequality corresponding to (\ref{ahm}). We note that if this inequality is strict then there exists a solution for the Yamabe problem on manifolds with boundary (see \cite{escobar:92}). Moreover, we show (in a forthcoming paper \cite{sr5}) that it also implies the existence of a spinor field solution of a Yamabe type boundary problem for the Dirac operator under the chiral bag boundary condition. A natural question is then to ask if one can find sufficient conditions for which this inequality is strict. It appears that the study of the Green function of the Dirac operator under this boundary condition allows to give such a condition.\\

Secondly, we will study the Green function for the Dirac operator for an other local elliptic boundary condition: the $\MIT$ bag boundary condition. Unlike the chiral bag boundary condition, this condition exists on every spin manifold with boundary since it does not require some additional structures. As in the previous part, we will compute the development of this function around a boundary point where the geometry of the manifold is quite simple. Then as an application, we will give a proof of the positive mass theorem proved by Escobar in \cite{escobar:92} which appears in the resolution of the Yamabe problem on manifolds with boundary. Note that our result requires that the manifold is spin (which is a stronger assumption compared with Escobar's theorem) however we don't assume that the whole manifold is locally conformally flat with umbilic boundary but only that there exists a boundary point $q\in\pa\M$ with a locally conformally flat with umbilic boundary neighbourhood (see below for the definition of such a neighbourhood).\\

{\bf Convention:} In this article, a point $q\in\pa\M$ in a Riemannian manifold $(\M^n,g)$ has a locally conformally flat with umbilic boundary neighbourhood, if there exist a metric $\ov{g}\in [g]$ and a neighbourhood $\V\subset\M$  (resp. $\U\subset\mathbb{R}^n_+$) of $q\in\pa\M$ (resp. of $0\in\pa\mathbb{R}^n_+$) such that $(\U,\xi)$ and $(\V,\ov{g})$ are isometric. In a same way, a point $q\in\pa\M$ has a locally flat with umbilic boundary neighbourhood if $\ov{g}=g$ in the previous notation.\\

{\bf Acknowledgements:} I would like to thank Oussama Hijazi and Emmanuel Humbert for their support. I am also very grateful to Marc Herzlich and Sebasti\'an Montiel for their remarks and their suggestions.


\section{Green function for the Dirac operator under the chiral bag boundary condition}\label{GFDOCBC}



\subsection{Definitions and first properties}\label{definition1}


In this section, we give a rigorous definition of the Green function of the Dirac operator under the chiral bag boundary condition. This condition has been introduced by Gibbons, Hawking, Horowitz and Perry \cite{hawking} (see also \cite{herzlich1}) on asymptotically flat manifolds with inner boundary in order to prove positive mass theorems for black holes. This condition defines an elliptic boundary condition for the Dirac operator and in the context of compact Riemannian spin manifolds with boundary, we can study the properties of its spectrum (see \cite{fs}, \cite{hmz} or \cite{sr} for example).\\
 
Let $(\M^n,g,\sigma)$ a $n$-dimensionnal connected compact Riemannian spin manifold with non empty smooth boundary $\pa\M$. We will denote by $\Sigma_g(\M)$ the spinor bundle over $(\M,g)$, $\na$ the Riemannian or spinorial Levi-Civita connections and ``$\cdot$'' the Clifford multiplication. The Dirac operator is then the first order elliptic differential operator acting on $\Sigma_g(\M)$ locally given by:
\begin{eqnarray*}
\D_g\varphi=\sum_{i=1}^n e_i\cdot\na_{e_i}\varphi,
\end{eqnarray*}
 
for all $\varphi\in\Ga\big(\Sigma_g(\M)\big)$ and where $\{e_1,...,e_n\}$ is a local $g$-orthonormal frame of the tangent bundle. We now briefly recall the definition of the chiral bag boundary condition. From now on, we assume that there exists a chirality operator, i.e. a linear map 
\begin{eqnarray*}
\Ga:\Sigma_g(\M)\longrightarrow\Sigma_g(\M)
\end{eqnarray*}

from the spinor bundle over $(\M,g)$ which satisfies the following properties:
\begin{equation}\label{poc}
\begin{array}{rc}
\Ga^2=\Id, & \<\Ga\varphi,\Ga\psi\>=\<\varphi,\psi\>\\
\na_X(\Ga\psi)=\Ga(\na_X\psi), & X\cdot\Ga\psi=-\Ga(X\cdot\psi),
\end{array}
\end{equation}

for all $X\in\Gamma(\TM)$ and for all spinor fields $\psi,\varphi\in\Gamma\big(\Sigma_g(\M)\big)$. If $\nu$ denotes the (inner) unit normal vector field to the boundary, an easy computation shows that the fiber preserving endomorphism:
\begin{eqnarray*}
v\cdot\Ga:\mathbf{S}_g\longrightarrow\mathbf{S}_g
\end{eqnarray*}

from the restricted spinor bundle $\mathbf{S}_g:=\Sigma_g(\M)_{|\pa\M}$ is an involution. So this spinor bundle splits into the direct sum $\mathbf{S}_g=\mathrm{V}_g^+\oplus\mathrm{V}_g^-$ where $\mathrm{V}_g^\pm$ is the eigensubbundle associated with the eigenvalues $\pm 1$. Thus one can check that the orthogonal projection 
$$\begin{array}{lccl}
\BCHIpm : & \L^2(\mathbf{S}_g) & \longrightarrow & \L^2(\V^{\pm}_g)\\ 
 & \varphi & \longmapsto & \frac{1}{2}(\Id\pm \nu\cdot\Ga)\varphi,
\end{array}$$ 

onto the eigensubbundle $\V^\pm_g$ defines an elliptic boundary condition for the Dirac operator. More precisely, the Dirac operator
\begin{eqnarray}
\D_g:\hpm_g=\{\varphi\in\H_1^2\;/ \;\BCHIpm(\varphi_{|\pa\M})=0\}\longrightarrow\L^2\big(\Si_g(\M)\big)
\end{eqnarray}

is a Fredholm operator and its spectrum (under this boundary condition) consists of entirely isolated real eigenvalues with finite multiplicity. In the following, we will denote by $\big(\lambda^\pm_k(g)\big)_{k\in\mathbb{Z}}$ this set of eigenvalues, i.e. we have:
$$\left\lbrace
\begin{array}{ll}
\D_g\varphi^\pm_k=\lambda_k^\pm(g)\varphi_k^\pm & \quad\textrm{on}\;\M\\ \\
\BCHIpm(\varphi^\pm_{k\,|\pa\M})=0 & \quad\textrm{along}\;\pa\M
\end{array}
\right.$$

where $(\phi^\pm_k)_{k\in\mathbb{Z}}$ can be chosen as being a spectral resolution. Let $\pi_1$, $\pi_2:\M\times\M\rightarrow\M$ be the projection on the first and the second component and let:
\begin{eqnarray*}
\Sigma_g(\M)\boxtimes\big(\Sigma_g(\M)\big)^\ast:=\pi_1^\ast\big(\Sigma_g(\M)\big)\otimes\big(\pi_2^*\big(\Sigma_g(\M)\big)\big)^\ast
\end{eqnarray*}

i.e. the fiber bundle whose fiber over $(x,y)$ is given by
\begin{eqnarray*} 
\Big(\Sigma_g(\M)\boxtimes\big(\Sigma_g(\M)\big)^\ast\Big)_{(x,y)}:=\mathrm{Hom}\big(\Sigma_{y}(\M),\Sigma_{x}(\M)\big),
\end{eqnarray*}

where $\Sigma_{y}(\M)$ denotes the fiber of the spinor bundle over $(\M,g)$. We are now ready to give the definition of the Green function for the Dirac operator under the chiral bag boundary condition.
\begin{definition}
A Green function for the Dirac operator under the chiral bag boundary condition (or a chiral Green function) is given by a smooth section:
\begin{eqnarray}
\G^\pm_\CHI:\M\times\M\setminus\Delta\rightarrow\Sigma_g(\M)\boxtimes\big(\Sigma_g(\M)\big)^\ast
\end{eqnarray}

locally integrable on $\M\times\M$ which satisfies, in a weak sense, the following boundary problem:
$$\left\lbrace
\begin{array}{ll}
\D_g\big(\G^\pm_\CHI(x,y)\big)=\delta_y\mathrm{Id}_{\Sigma_y\M}\\ \\
\BCHI^\pm\big(\G^\pm_\CHI(x,y)\big)=0,\quad\textrm{for}\;x\in\pa\M\setminus\{y\},
\end{array}
\right.$$

for all $x\neq y\in\M$ (where $\Delta:=\{(x,y)\in\M\times\M\,/\,x=y\}$). In other words, we have:
\begin{eqnarray*}
\int_{\M}\<\G_{\CHI}^\pm(x,y)\psi_0^\pm,\D_g\varphi(x)\>dv(x)=\<\psi_0^\pm,\varphi(y)\>
\end{eqnarray*}

for all $y\in\M$, $\psi_0^\pm\in\Sigma_{y}(\M)$ satisfying $\nu\cdot\Ga\psi_0^\pm=\mp\psi_0^\pm$ and $\varphi\in\Gamma\big(\Sigma_g(\M)\big)$ such that $\BCHIpm(\varphi_{|\pa\M})=0$.
\end{definition}

\begin{remark}\label{remark1}
Here we recall that on the Euclidean space, we have a corresponding notion of Green function for the Dirac operator (see \cite{amm4} for example). In fact, one can easily check that there exists a unique Green function given by:
$$\begin{array}{lccl}
\mathrm{G}_{\eucl}: & \mathbb{R}^n\times\mathbb{R}^n\setminus\Delta & \longrightarrow & \Sigma\mathbb{R}^n\boxtimes\big(\Sigma\mathbb{R}^n\big)^\ast\\
 & (x,y) & \longmapsto & -\frac{1}{\omega_{n-1}}\frac{x-y}{|x-y|^n}\cdot,
\end{array}$$ 

where $\omega_n$ stands for the volume of the $n$-dimensional round sphere. This Green function can be seen (up to a conformal change of metric) as the inverse of the Dirac operator on the standard sphere $\mathbb{S}^n$.
\end{remark}

We are now going to give the development of the chiral Green function near a boundary point $q\in\pa\M$ and when the manifold has a particular geometry near this point. From now on, we assume that the Dirac operator is invertible under the chiral bag boundary condition, i.e. $\mathrm{Ker}^\pm(\D_g)=\{0\}$ where $\mathrm{Ker}^\pm(\D_g)$ denotes the kernel of the Dirac operator with domain in $\mathcal{H}^\pm_g$.\\

In order to get this expansion, we have to identify (locally) the spinor bundle over $(\M,g)$ with the (trivial) spinor bundle over an open set of the half-space $\mathbb{R}^n_+$ endowed with the Euclidean metric. In fact, in a previous work (see \cite{sr3} or \cite{mathese}), it is shown that if 
\begin{eqnarray*}
\mathcal{F}_q:\U\longrightarrow\V
\end{eqnarray*}

stands for the Fermi coordinate system around $q\in\V$ where $\U$ (resp. $\V$) is an open set in $\mathbb{R}^n_+$ (resp. in $\M$) then we have a canonical trivialization given by:
$$\begin{array}{clc}
\Sigma_\xi(\U) & \longrightarrow & \Sigma_g(\V)  \\
 \psi & \longmapsto & \overline{\psi}.
\end{array}$$

where $\xi$ is the euclidean metric on $\mathbb{R}^n_+$. This identification is closely related to those given in \cite{bourguignon} and \cite{amm3}. 

\begin{remark}\label{remark2}
In \cite{sr3} (see also \cite{mathese}), we show that it is sufficient to choose a constant spinor field $\psi_0\in\Sigma_0(\mathbb{R}^n_+)$ such that $\BCHIpm\big(\ov{\psi_0^\pm}\big)(q)=0$ at one point $q\in\pa\M$ so that $\BCHIpm\big(\ov{\psi_{0}^\pm}_{\,|\V\cap\pa\M}\big)=0$
\end{remark}

In the following, we will note indifferently a spinor over $\U$ and its image on $\V$ under this trivialization. The expansion of the chiral Green function is then given by:
\begin{proposition}\label{Devfoncgreen}
Assume that there exists a point $q\in\pa\M$ which has a locally flat with umbilic boundary neighbourhood. Then the chiral Green function exists and has, in the above trivialization, the following expansion (near $q$):
\begin{eqnarray}\label{devfoncgreen}
\G_\CHI^\pm(x,q)\psi_0^\pm=-\frac{2}{\omega_{n-1}}\frac{x-q}{|x-q|^n}\cdot\psi_0^\pm+\mathrm{m}_\CHI^\pm(x,q)\psi_0^\pm,
\end{eqnarray}

where $\psi_0^\pm\in\Sigma_q(\M)$ is a spinor field satisfying $\nu\cdot\Ga\psi_0^\pm=\mp\psi_0^\pm$ (in other words $\BCHI^{\pm}(\psi_0^\pm)=0$) and where $\mathrm{m}_\CHI^\pm(\,.\,,q)\psi_0^\pm$ is a smooth spinor field such that $\D_g\big(\mathrm{m}_\CHI^\pm(\,.\,,q)\psi_0^\pm\big)=0$ around $q$. Moreover, along $\pa\M$, we have:
\begin{eqnarray}\label{masschiral}
\BCHI^{\pm}\big(\mathrm{m}_\CHI^\pm(.\,,q)\psi^\pm_{0|\pa\M}\big)=0.
\end{eqnarray}
\end{proposition}

{\it Proof:} 
We only prove this result for the boundary condition $\mathbb{B}_{\CHI}^-$ since the proof is similar for $\mathbb{B}_{\CHI}^+$. We trivialize the spinor bundle over $\V$ (around the boundary point $q$) using the preceding discussion. In fact, we can isometrically identify $(\V,g)$ with $(\U,\xi)$ and thus the spinor bundles over these spaces. Now we consider a smooth cut-off function $\eta$ such that:
$$\left\lbrace
\begin{array}{ll}
0\leq \eta\leq 1, & \eta\in\mathrm{C}^\infty(\M)\\
 \eta\equiv 1\textrm{ on } \mathrm{B}^+_q(\delta), & \mathrm{supp}(\eta)\subset\mathrm{B}^+_q(2\delta)\subset\V
\end{array}
\right.$$

where $\mathrm{B}^+_q(2\delta)$ is the half-ball centered in $q$ and with radius $2\delta>0$ contained in $\V$. Let $\Psi=\D_g\big(\eta\G_{\eucl}(\,.\,,q)\psi_0\big)$ with $\psi_0\in\Sigma_{q}(\M)$ such that $\nu\cdot\Ga\psi_0=\psi_0$. Using Remark~\ref{remark2}, we then have that this boundary condition is fulfilled along $\V\cap\pa\M$ for the spinor field $\psi_0$. The spinor field $\Psi$ is smooth on $\M\setminus\{q\}$ and an easy calculation leads to $\Psi_{|\mathrm{B}_q^+(\delta)\setminus\{q\}}=0$. Thus we can extend this spinor all over $\M$ by $\Psi(q)=0$. Since the Dirac operator is supposed to be invertible on $\mathcal{H}^-$, there exists a unique spinor field $\mathrm{m}_\CHI^-(\,.\,,q)\psi_0\in\mathcal{H}^-$ such that $\D_g\big(\mathrm{m}_\CHI^-(\,.\,,q)\psi_0\big)=-\Psi$. Now for $x\in\M\setminus\{q\}$, we let:
\begin{eqnarray*}
\G^-_q(x)\psi_0=2\eta\,\G_{\eucl}(x,q)\psi_0+\mathrm{m}_\CHI^-(x,q)\psi_0.
\end{eqnarray*}

We first check that if $x\in\pa\M\setminus\{q\}$ then this spinor field satisfies the chiral bag boundary condition. If $x\in\mathrm{supp}(\eta)^c\cap\pa\M$ then
\begin{eqnarray*}
\BCHI^{-}\big(\G^-_q(x)\psi_0\big)=\BCHI^{-}\big(\mathrm{m}_\CHI^-(x,q)\psi_0\big)=0
\end{eqnarray*}

since by construction $\mathrm{m}_\CHI^\pm(\,.\,,q)\psi_0\in\mathcal{H}^-$. Secondly, if $x\in\big(\mathrm{supp}(\eta)\cap\pa\M\big)\setminus\{q\}$ we get
\begin{eqnarray*}
\BCHI^{-}\big(\G^-_q(x)\psi_0\big)=\BCHI^{-}\big(\G_{\eucl}(x,q)\psi_0\big).
\end{eqnarray*}

Using the expression of the Euclidean Green function given in Remark~\ref{remark1} and the properties (\ref{poc}) of the chirality operator $\Gamma$, we obtain
\begin{eqnarray*}
\nu\cdot\Ga\big(\frac{x-q}{|x-q|^n}\cdot\psi_0\big)=\frac{x-q}{|x-q|^n}\cdot\nu\cdot\Ga\psi_0=\frac{x-q}{|x-q|^n}\cdot\psi_0
\end{eqnarray*}

since $\psi_0$ is chosen such that $\nu\cdot\Ga\psi_0=\psi_0$ and thus $\mathbb{B}_{\CHI}^-\big(\G^-_q(x)\psi_0\big)=0$. An easy calculation shows that for all $\phi\in\mathcal{H}^-$, we have
\begin{eqnarray*}
\int_{\M}\<\G^-_q(x)\psi_0,\D(\phi)\>dv(x)=\<\psi_0,\phi(q)\>,
\end{eqnarray*}

and so $\G^-_q(\,.\,)$ is a chiral Green function. Unicity follows from the hypothesis $\mathrm{ker}^\pm(\D_g)=\{0\}$. In fact, if we assume that there exists another chiral Green function $\widetilde{\G}^-_q(\,.\,)$, we can show that the spinor field $\G^-_q(\,.\,)\psi_0-\widetilde{\G}^-_q(\,.\,)\psi_0$ is in the kernel of the Dirac operator (using the classical regularity theorems, see \cite{schwartz} for example). Since the Dirac operator (under the chiral bag boundary condition) is supposed to be invertible, this spinor vanishes identically and so unicity follows directly.
\hfill$\square$\\

In the following, the chiral Green function will be indifferently denoted by $\G_\CHI^\pm(\,.\,,q)$ or $\G^\pm_q(\,.\,)$.\\

We now look at the behaviour of the chiral Green function under a conformal change of metric. In fact, we prove the following result:
\begin{proposition}\label{Confgreen}
Let $\ov{g}=f^2g$ be a metric in the conformal class of $g$. If $\G_\CHI^\pm(\,.\,,q)$ (resp. $\ov{\G}_\CHI^\pm$) denote the chiral Green function for the Dirac operator $\D_g$ (resp. $\D_{\ov{g}}$), then:
\begin{eqnarray}\label{covgreen}
\ov{\G}_\CHI^\pm(\,.\,,q)=f^{-\frac{n-1}{2}}(\,.\,)\,f^{-\frac{n-1}{2}}(q)\,\ov{\G_\CHI^{\pm}}(\,.\,,q)
\end{eqnarray}

for $q\in\pa\M$.
\end{proposition}

{\it Proof:}
First recall that given two conformal metrics $g$ and $\ov{g}=f^2g$, there exists a canonical identification between the spinor bundle over $(\M,g)$ and the one over $(\M,\ov{g})$ (see \cite{hitchin:74} or \cite{hijazi:86}). It is given by the bundle isomorphism (which is a fiberwise isometry)
\begin{eqnarray}
\F:\Sigma_g(\M)\longrightarrow\Sigma_{\ov{g}}(\M)
\end{eqnarray}

such that the Dirac operators satisfy the relation
\begin{eqnarray}\label{covconf}
\F^{-1}\big(\D_{\ov{g}}(\F(\psi))\big)=f^{-\frac{n+1}{2}}\D_g(f^{\frac{n-1}{2}}\psi),
\end{eqnarray}

for all $\psi\in\Gamma\big(\Sigma_g(\M)\big)$. On the other hand, since $\G_q^{-}$ is the chiral Green function for the Dirac operator, we have:
\begin{eqnarray*}
\int_{\M}\<\G^-_q(x)\psi_0,\D_g\varphi\>dv(g)=\<\psi_0,\varphi(q)\>
\end{eqnarray*}

for all $\psi_0\in\Sigma_q(\M)$ such that $\BCHI^{-}(\psi_0)=0$ and $\varphi\in\Ga\big(\Sigma_g(\M)\big)$ such that $\BCHI^{-}(\varphi_{|\pa\M})=0$. We can now prove that the section defined by (\ref{covgreen}) is a chiral Green function for the Dirac operator $\D_{\ov{g}}$. For $\varphi\in\Gamma\big(\Sigma_g(\M)\big)$ such that $\BCHI^{-}(\varphi_{|\pa\M})=0$, we let $\Phi=f^{-\frac{n-1}{2}}\F(\varphi)\in\Gamma\big(\Sigma_{\ov{g}}(\M)\big)$. This spinor field satisfies $\ov{\mathbb{B}}_\CHI^\pm(\phi_{|\pa\M})=0$ because of the conformal covariance of the chiral bag boundary condition (see for example \cite{mathese}). We now write:
\begin{eqnarray*}
\int_{\M}\<\ov{\G}^{-}_q(x)\F(\psi_0),\D_{\ov{g}}\Phi\>dv(\ov{g}) & = & \int_{\M}f^{-\frac{n-1}{2}}(x)\,f^{-\frac{n-1}{2}}(q)f^{-\frac{n+1}{2}}(x)\<\G^{-}_q(x)\psi_0,\D_g\varphi\>f^n dv(g)\\
& = & f^{-\frac{n-1}{2}}(q)\int_{\M}\<\G^{-}_q(x)\psi_0,\D_g\varphi\>dv(g)\\
& = & \<\psi_0,f^{-\frac{n-1}{2}}(q)\varphi(q)\>
\end{eqnarray*}

that is:
\begin{eqnarray*}
\int_{\M}\<\ov{\G}^{-}_q(x)\F(\psi_0),\D_{\ov{g}}(\Phi)\>dv(\ov{g}) & = & \<\F(\psi_0),\Phi(q)\>
\end{eqnarray*}

for all $\F(\psi_0^\pm)\in\Sigma_q(\M)$ such that $\ov{\mathbb{B}}_\CHI^\pm\big(\F(\psi_0^\pm)\big)=0$ and $\phi\in\Gamma\big(\Sigma_{\ov{g}}(\M)\big)$ such that $\ov{\mathbb{B}}_\CHI^\pm(\phi_{|\pa\M})=0$. We have then checked that $\ov{\G}^{-}_q(\,,\,)$ is a chiral Green function. Unicity follows directly.
\hfill$\square$

\begin{remark}
Propositions \ref{Devfoncgreen} and \ref{Confgreen} imply that the chiral Green function exists on locally conformally flat manifolds with umbilic boundary on which the Dirac operator is invertible under the chiral bag boundary condition.
\end{remark}

We now prove a self-adjointness result for the chiral bag boundary condition which will be very important for the following and more particulary for the application given in Section \ref{apl1}.
\begin{proposition}\label{greensymetry}
For $(x,y)\in\M\times\M\setminus\Delta$, we have:
\begin{eqnarray}\label{greensym}
\big(\G_\CHI^{\pm}\big)^*(x,y)= \G_\CHI^{\pm}(y,x).
\end{eqnarray}

In fact, if $\psi_0^\pm\in\Sigma_x(\M)$ (resp. $\varphi_0^\pm\in\Sigma_y(\M)$) such that $\nu\cdot\Ga\psi_0^\pm=\mp\psi_0^\pm$ for $x\in\pa\M$ (resp. $\nu\cdot\Ga\varphi_0^\pm=\mp\varphi_0^\pm$ for $y\in\pa\M$), then:
\begin{eqnarray*}
\<\G_\CHI^{\pm}(x,y)\varphi_0^\pm,\psi_0^\pm\>_{\Sigma_x(\M)}=\<\varphi_0^\pm,\G^{\pm}_\CHI(y,x)\psi_0^\pm\>_{\Sigma_y(\M)}.
\end{eqnarray*}
\end{proposition}

{\it Proof:}
We use the fact that there exists a spectral resolution of the space of $\L^2$-spinors, i.e. for all $\psi\in\Gamma\big(\Sigma_g(\M)\big)$, $\exists (\mathrm{A}_k)_{k\in\mathbb{Z}}\subset\mathbb{C}$ such that $\psi=\sum_{k\in\mathbb{Z}}\mathrm{A}_k\varphi_k$ where $\varphi_k$ is a smooth spinor field satisfying
$$\left\lbrace
\begin{array}{ll}
\D_g\varphi_k=\lambda_k(g)\varphi_k & \quad\textrm{on}\;\M\\
\BCHI^{-}(\varphi_{k\,|\pa\M})=0 & \quad\textrm{along}\;\pa\M.
\end{array}
\right.$$

First note that:
 \begin{eqnarray*}
\int_{\M\times\M\setminus\Delta}\<\G_\CHI^-(x,y)\varphi_i,\varphi_j\>dv(x)dv(y) & = & \int_{\M}\Big(\int_{\M}\<\G_\CHI^-(x,y)\varphi_i,\varphi_j\>dv(x)\Big)dv(y)\\
& = & \int_{\M}\Big(\int_{\M}\<\G_\CHI^-(x,y)\varphi_i,\frac{1}{\lambda_j}\D_g\varphi_j\>dv(x)\Big)dv(y)\\
& = & \frac{1}{\lambda_j}\int_{\M}\<\varphi_i,\varphi_j\>dv(y)\\
& = & \frac{1}{\lambda_j}\delta_{ij},
\end{eqnarray*}

since $\D$ is assumed to be invertible and thus $\lambda_j\neq 0$ for all $j$. In a same way, we easily show that:
\begin{eqnarray*}
\int_{\M\times\M\setminus\Delta}\<\varphi_i,\G_\CHI^-(y,x)\varphi_j\>dv(x)dv(y) & = & \frac{1}{\lambda_i}\delta_{ij},
\end{eqnarray*}

We thus have proved that for all $(i,j)\in\mathbb{Z}^2$:
\begin{eqnarray*}
\int_{\M\times\M\setminus\Delta}\<\G_\CHI^-(x,y)\varphi_i,\varphi_j\>dv(x)dv(y)=\int_{\M\times\M\setminus\Delta}\<\varphi_i,\G_\CHI^-(y,x)\varphi_j\>dv(x)dv(y).
\end{eqnarray*}

By linearity, we have:
\begin{eqnarray*}
\int_{\M\times\M\setminus\Delta}\<\G_\CHI^-(x,y)\varphi,\psi\>dv(x)dv(y)=\int_{\M\times\M\setminus\Delta}\<\varphi,\G_\CHI^-(y,x)\psi\>dv(x)dv(y),
\end{eqnarray*}

for all $\psi$, $\varphi\in\Gamma\big(\Sigma_g(\M)\big))$ such that $\BCHI^-(\psi_{|\pa\M})=\BCHI^-(\varphi_{|\pa\M})=0$. Now we consider $x_0\in\M$ such that $y\mapsto\G_\CHI^-(y,x_0)$ and let $f_{\varepsilon}:\mathbb{R}^n\rightarrow\mathbb{R}$ such that $f_\varepsilon\rightarrow\delta_{x_0}$ where $\delta_{x_0}$ is the Dirac's distribution at $x_0$. So if $\psi\in\Gamma\big(\Sigma_g(\M)\big)$ satisfies $\BCHI^-(\psi_{|\pa\M})=0$ and $\psi(x_0)=\psi_0$, we obtain:
\begin{eqnarray*}
\int_{\M}\<\varphi(y),\G_\CHI^-(y,x_0)\psi_0\>dv(y)=\int_{\M}\<\G_\CHI^-(x_0,y)\varphi(y),\psi_0\>dv(y).
\end{eqnarray*}

Using the same method for the $y-$variable with $y_0\neq x_0$ leads to:
\begin{eqnarray*}
\<\varphi_0,\G_\CHI^-(y_0,x_0)\psi_0\>_{\Sigma_{y_{0}}(\M)}=\<\G_\CHI^-(x_0,y_0)\varphi_0,\psi_0\>_{\Sigma_{x_{0}}(\M)},
\end{eqnarray*}

and this concludes the proof.
\hfill$\square$\\

We now show that the chiral Green functions $\G_\CHI^-(\,.\,,q)$ and $\G_\CHI^+(\,.\,,q)$ are related by the action of the chirality operator. In fact, we prove the following result:
\begin{proposition}\label{posneg}
Let $(\M^n,g)$ be a connected compact Riemannian spin manifold with a boundary point $q\in\pa\M$ which has a locally conformally flat with umbilic boundary neighbourhood. Then we get:
\begin{eqnarray}
\G_\CHI^-(\,.\,,q)=-\Ga\G_\CHI^+(\,.\,,q)\Ga.
\end{eqnarray}
\end{proposition}

{\it Proof:}
Without loss of generality, we can assume that the metric $g$ is such that there is a neighbourhood of $q\in\pa\M$ isometric to an open set of the Euclidean half-space. Then in this chart, the chiral Green function $\G_\CHI^-(\,.\,,q)$ admits the development (\ref{devfoncgreen}). Now define the section $\widetilde{\G}(\,.\,,q)=-\Ga\G_\CHI^+(\,.\,,q)\Ga$ on $\M\setminus\{q\}$. First note that we can easily check that $\widetilde{\G}(\,.\,,q)\in\mathcal{H}^-$. Then a direct computation shows that the spinor field given by
\begin{eqnarray*}
\G_\CHI^-(\,.\,,q)\psi_0-\widetilde{\G}(\,.\,,q)\psi_0 
\end{eqnarray*}

for $\psi_0$ such that $\BCHI^-(\psi_0)=0$, is harmonic. Hence since the Dirac operator is supposed to be invertible under the chiral bag boundary condition, we conclude that $\widetilde{\G}(\,.\,,q)\psi_0=\G_\CHI^-(\,.\,,q)\psi_0$ is the unique chiral Green function. 
\hfill$\square$\\

We are now ready to define the chiral mass operator (see \cite{amm4} for the case of manifolds without boundary). The name of this operator comes from its tight relation with the concept of mass in General Relativity which also appears in the context of the Yamabe problem (see \cite{schoen:84} or \cite{lee.parker:87} for the closed case and \cite{escobar:92} for the boundary case). 

\begin{definition}
Let $(\M^n,g,\sigma)$ be a connected compact Riemannian spin manifold with non empty smooth boundary $\pa\M$. Suppose that there exits a boundary point $q\in\pa\M$ which has a locally flat with umbilic boundary neighbourhood. The chiral mass operator is then defined by:
$$\begin{array}{rll}
\mathrm{m}_\CHI^\pm(q):\V^\pm_q & \longrightarrow & \V^\pm_q \\
\psi_0^\pm & \longmapsto & \mathrm{m}^\pm_\CHI(q,q)\psi_0^\pm
\end{array}$$

where $\mathrm{m}^\pm(\,.\,,q)\psi_0^\pm$ is given in Proposition~\ref{Devfoncgreen} and $\V^\pm_q=\frac{1}{2}\big(\mathrm{Id}\pm\nu\cdot\Ga\big)\big(\Sigma_q(\M)\big)$. 
\end{definition}

Thanks to Proposition \ref{greensym}, we can deduce the following result which will be very useful for the next section.
\begin{proposition}
For a point $q\in\pa\M$ with a locally flat with umbilic boundary neighbourhood, the chiral mass operator $\mathrm{m}_\CHI^\pm(q)$ is linear and symmetric.
\end{proposition}

{\it Proof:}
An easy computation shows that the chiral Green function is linear. Using this fact, it clearly follows that the chiral mass operator is also linear. The symmetry of the chiral mass operator comes from the symmetry of the chiral Green function proved in Proposition~\ref{greensymetry}.
\hfill$\square$\\

\begin{remark}
In the following, we will refer to the ``negative chiral mass operator'' (resp. ``positive chiral mass operator'') for $\mathrm{m}_\CHI^-(q)$ (resp. $\mathrm{m}_\CHI^+(q)$).
\end{remark}

As a direct consequence, we obtain:
\begin{corollary}\label{spectrum}
For a point $q\in\pa\M$ with a locally flat with umbilic boundary neighbourhood, the (pointwise) spectrum of the chiral mass operator $\mathrm{m}_\CHI^\pm(q)$ is real. Moreover, if $\kappa$ is an eigenvalue for the negative chiral mass operator, then $-\kappa$ is an eigenvalue for the positive chiral mass operator.
\end{corollary}

{\it Proof:}
Using Proposition~\ref{posneg}, we can easily check that the positive and the negative chiral mass operators satisfy the relation:
\begin{eqnarray*}
\Ga\mathrm{m}_\CHI^-(q)=-\mathrm{m}_\CHI^+(q)\Ga.
\end{eqnarray*}

So consider an eigenspinor $\psi_0\in\V^-_q$ for the negative chiral mass operator associated with the eigenvalue $\kappa$, i.e. $\mathrm{m}_\CHI^-(q)\psi_0=\kappa\psi_0$ and $\nu\cdot\Ga\psi_0=\psi_0$. Using the preceding formula and the fact that $\nu\cdot\Ga\big(\Ga\psi_0\big)=-\Ga\psi_0$, we observe that $\mathrm{m}_\CHI^+(q)(\Ga\psi_0)=-\kappa\Ga\psi_0$ and thus $-\kappa$ is an eigenvalue for the positive chiral mass operator.
\hfill$\square$\\

In the next section, we give a direct application of the construction of the chiral Green function and the chiral mass operator.


\subsection{Application: The chiral bag invariant}\label{apl1}


This part is devoted to a direct application of the preceding construction of the chiral Green function for the Dirac operator. This application concerns a spin conformal invariant on manifolds with boundary introduced in \cite{sr3}. We first begin with a brief introduction on this invariant. We have seen in Section~\ref{definition1} (see \cite{hijazi.montiel.roldan:01} for more details) that the spectrum of the Dirac operator under the chiral bag boundary condition consists of entirely isolated real eigenvalues with finite multiplicity. If we denote by $\lambda_1^\pm(g)$ the first eigenvalue of the Dirac operator $\D_g$ under the boundary condition $\BCHI^\pm$, then the chiral bag invariant is defined by:
\begin{eqnarray}\label{cbi}
\lambda_{\min}^{\pm}(\M,\pa\M):=\underset{\ov{g}\in[g]}{\inf}|\lambda^{\pm}_1(\ov{g})|\Vol(\M,\ov{g})^{\frac{1}{n}},
\end{eqnarray}

where $[g]$ denotes the conformal class of $g$ and $\Vol(\M,\ov{g})$ is the volume of the manifold $\M$ equipped with the Riemannian metric $\ov{g}\in[g]$. A very useful formula for the following is given by the variational characterization of the chiral bag invariant. In fact, it is shown in \cite{sr3} (see also \cite{ammann}) that:
\begin{eqnarray}\label{charvar}
\lambda_{\min}^\pm(\M,\pa\M)=
\underset{\varphi\in\Coker^\pm_g}{\inf}\Big\{\frac{\big(\int_{\M}|\D_g\varphi|^{\frac{2n}{n+1}}dv(g)\big)^{\frac{n+1}{n}}}{\big|\int_{\M}\mathrm{Re}\<\D_g\varphi,\varphi\>dv(g)\big|}\Big\},
\end{eqnarray}

where $\Coker^\pm_g$ is the $\L^2$-orthogonal of $\Ker^{\pm}(\D_g)$ in $\mathcal{H}^\pm_g$ and $\D_g$ is the Dirac operator in the metric $g$.
\begin{remark}
\begin{enumerate}

\item The above definition seems to depend on the boundary condition chosen $\mathbb{B}^{+}_{g}$ or $\mathbb{B}^{-}_{g}$, however it doesn't (see \cite{sr3}). Taking in account this fact, we will denote by $\lambda_{\min}(\M,\pa\M)$ the chiral bag invariant in what follows and we will use the $\mathbb{B}^-_g$ condition. 

\item Using the Hijazi inequality proved in \cite{sr}, we can compare the chiral bag invariant with the Yamabe invariant $\mu(\M,\pa\M)$ of the manifold (see \cite{escobar:92}). In fact, if $n\geq 3$, we have:
\begin{eqnarray}\label{hijbord}
\lambda_{\mathrm{min}}(\M,\pa\M)^2\geq\frac{n}{4(n-1)}\mu(\M,\pa\M).
\end{eqnarray}

\item In \cite{sr3}, we have shown a spinorial analogous of Aubin's (or Escobar's) inequality \cite{aubin:76} (\cite{escobar:92}, for the non empty boundary case). More precisely, we proved that if $n\geq 2$:
\begin{eqnarray}\label{largeb1}
\lambda_{\min}(\M,\pa\M)\leq\lambda_{\min}(\hs,\pa\hs)=\frac{n}{2}\Big(\frac{\omega_n}{2}\Big)^{\frac{1}{n}}.
\end{eqnarray}
 
\noindent This inequality is the analogue of the one obtained by Ammann, Humbert and Morel \cite{amm3} in the boundaryless case.
\end{enumerate}
\end{remark}

A natural question is to find a sufficient condition under which Inequality~(\ref{largeb1}) is strict. Keeping in mind the work of Schoen \cite{schoen:84}, Escobar \cite{escobar:92} and Ammann, Humbert and Morel \cite{amm4}, we are going to see that the construction of the chiral mass operator gives an answer for a certain class of manifolds. More precisely, we prove the following result:
\begin{theorem}\label{lcfm}
Let $(\M^n,g,\sigma)$ be a $n$-dimensional $(n\geq 2)$ connected compact Riemannian spin manifold with non empty smooth boundary $\pa\M$. Suppose that there exists a point $q\in\pa\M$ which has a locally conformally flat with umbilic boundary neighbourhood. Moreover, we assume that the Dirac operator $\D_g$ is invertible on $\mathcal{H}^-$ (or on $\mathcal{H}^+$) and that the negative chiral mass operator $\mathrm{m}_\CHI^-(q)$ (or the positive chiral mass operator $\mathrm{m}_\CHI^+(q)$) is not identically zero. Then we get:
\begin{eqnarray*}
\lambda_{\min}(\M,\pa\M)<\lambda_{\min}(\hs,\pa\hs).
\end{eqnarray*}
\end{theorem}

The proof of this theorem is in the same spirit as the one of Inequality~(\ref{largeb1}). Indeed, the idea is to construct a test-spinor to estimate in the variational characterization~(\ref{charvar}) of $\lambda_{\min}(\M,\pa\M)$. In \cite{sr3}, the test-spinor was constructed from a Killing spinor on the hemisphere satisfying the chiral bag boundary condition and with support contained in an open set of a trivialization around a boundary point. In order to prove Theorem~\ref{lcfm}, it is not enough to extend by zero the test-spinor away from the open set of trivialization. In fact, the good extention is given by the chiral Green function. More precisely, we show:
\begin{theorem}\label{cop}
Let $(\M,g,\sigma)$ be a $n$-dimensional connected compact Riemannian spin manifold ($n\geq 2$) with non empty boundary $\pa\M$. Suppose that there exists a point $q\in\pa\M$ which has a locally conformally flat with umbilic boundary neighbourhood. Assume also that there exits on $\M\setminus\{q\}$ a spinor field $\psi^\pm$ such that:
\begin{equation}\label{hyp1}
\left\lbrace
\begin{array}{ll}
\D_g\psi^\pm=0 & \quad\text{on}\;\M\setminus\{q\}\\ \\
\BCHI^\pm(\psi^\pm_{|\pa\M})=0 & \quad\text{along}\;\pa\M
\end{array}
\right.
\end{equation}

which admits the following development around $q$:
\begin{eqnarray}\label{hyp2}
\psi^\pm=\frac{x}{r^n}\cdot\psi_0^\pm+\psi_1^\pm+\theta^\pm,
\end{eqnarray}

where $\psi_0^\pm$, $\psi_1^\pm\in\Sigma_q\M$ are spinors satisfying $\nu\cdot\Ga\psi_0^\pm=\mp\psi_0^\pm$, $\nu\cdot\Ga\psi_1^\pm=\mp\psi_1^\pm$ and such that:
\begin{eqnarray}\label{hyp3}
\mathrm{Re}\,\<\psi_0^\pm,\psi_1^\pm\><0\qquad\text{and}\qquad\mathrm{Re}\,\<x\cdot\psi_0^\pm,\psi_1^\pm\>=0.
\end{eqnarray}

We also assume that $\theta^\pm$ is a smooth spinor field all over $\M$ satisfying $\theta^\pm=O(r)$, $\BCHI^\pm(\theta^\pm_{|\pa\M})=0$ and which is harmonic around $q$. Under these hypothesis, we get:
\begin{eqnarray*}
\lambda_{\min}(\M,\pa\M)<\lambda_{\min}(\hs,\pa\hs)=\frac{n}{2}\Big(\frac{\omega_n}{2}\Big)^{\frac{1}{n}}.
\end{eqnarray*}
\end{theorem}

{\it Proof:}
For $\varepsilon>0$, we let $\rho:=\varepsilon^{\frac{1}{n+1}}$ and $\varepsilon_0:=\frac{\rho^n}{\varepsilon}f\big(\frac{x}{\varepsilon}\big)^{\frac{n}{2}}$ with $f(r)=\frac{1}{1+r^2}$ and we consider the spinor field defined by:
\begin{equation}
\psi_{\varepsilon}^\pm=
\left\lbrace
\begin{array}{ll}
f\big(\frac{x}{\varepsilon}\big)\big(1-\frac{x}{\varepsilon}\big)\cdot\psi_0^\pm-\varepsilon_0\psi_1^\pm \quad & \quad \textrm{if }r\leq \rho\\ \\
 -\varepsilon_0\big(\psi^\pm-\eta\theta^\pm\big)+\eta f\big(\frac{\rho}{\varepsilon}\big)^{\frac{n}{2}}\psi_0^\pm\quad & \quad\textrm{if } \rho\leq r\leq 2\rho\\ \\
 \varepsilon_0\psi^\pm \quad & \quad\textrm{if }2\rho\leq r
\end{array}
\right.
\end{equation}

where $r=d(x,q)$ and $\eta$ is a cut-off function equal to zero on $\M\setminus\mathrm{B}^+_q(2\rho)$, $1$ on $\mathrm{B}^+_q(\rho)$ and satisfying $|\na\eta|\leq\frac{2}{\rho}$. We first check that $\BCHI^\pm(\psi^\pm_{\varepsilon|\pa\M})=0$. For $r\leq \rho$, we have:
\begin{eqnarray*}
\BCHI^\pm(\psi^\pm_{\varepsilon|\pa\M}) & = & \BCHI^\pm\Big(\big(f\big(\frac{x}{\varepsilon}\big)\big(1-\frac{x}{\varepsilon}\big)\cdot\psi_0^\pm-\varepsilon_0\psi^\pm_{1}\big)_{|\pa\M}\Big)\\ \\
& = & f\big(\frac{x}{\varepsilon}\big)\big(1-\frac{x}{\varepsilon}\big)\cdot\,\BCHI^\pm(\psi_{0\,|\pa\M}^\pm)-\varepsilon_0\,\BCHI^\pm(\psi_{1\,|\pa\M}^\pm)=0
\end{eqnarray*}

since $\BCHI^\pm(\psi_{0\,|\pa\M}^\pm)=\BCHI^\pm(\psi_{1\,|\pa\M}^\pm)=0$. In the same way, since $\psi^\pm$ and $\theta^\pm$ satisfy $\BCHI^\pm(\psi_{|\pa\M}^\pm)=\BCHI^\pm(\theta_{|\pa\M}^\pm)=0$, we easily check that for $\rho\leq r\leq 2\rho$ or $r\geq 2\rho$, we also have $\BCHI^\pm(\psi^\pm_{\varepsilon|\pa\M})=0$. Without loss of generality, we can assume that $|\psi_0^\pm|^2=1$. Since $\psi^\pm$ and $\theta^\pm$ are harmonic around $q$, we compute that:
\begin{equation}
\D_g\psi_{\varepsilon}^\pm=
\left\lbrace
\begin{array}{ll}
\frac{n}{\varepsilon}f\big(\frac{r}{\varepsilon}\big)^{\frac{n}{2}+1}\big(1-\frac{x}{\varepsilon}\big)\cdot\psi_0^\pm\quad & \quad \textrm{if }r\leq \rho\\ \\
 \varepsilon_0\na\eta\cdot\theta^\pm+f\big(\frac{r}{\varepsilon}\big)^{\frac{n}{2}}\na\eta\cdot\psi_0^\pm\quad & \quad\textrm{if } \rho\leq r\leq 2\rho\nonumber\\ \\
 0\quad & \quad\textrm{if }2\rho\leq r.
\end{array}
\right.
\end{equation}
 
In order to obtain an estimate of the numerator of the variational characterization~(\ref{charvar}) of $\lambda_{\min}(\M,\pa\M)$, one can check that:
\begin{equation}
|\D_g\psi^\pm_{\varepsilon}|^{\frac{2n}{n+1}}=
\left\lbrace
\begin{array}{ll}
n^{\frac{2n}{n+1}}\varepsilon^{-\frac{2n}{n+1}}f\big(\frac{r}{\varepsilon}\big)^{n}\quad & \quad \textrm{if }r\leq \rho\\ \\
 |\varepsilon_0\na\eta\cdot\theta^\pm+f\big(\frac{r}{\varepsilon}\big)^{\frac{n}{2}}\na\eta\cdot\psi_0^\pm|^{\frac{2n}{n+1}}\quad & \quad\textrm{if } \rho\leq r\leq 2\rho\nonumber\\ \\
 0\quad & \quad\textrm{if }2\rho\leq r.
\end{array}
\right.
\end{equation}

Hence we have:
\begin{eqnarray*}
\int_{\mathrm{B}^+_q(\rho)}|\D_g\psi^\pm_{\varepsilon}|^{\frac{2n}{n+1}}dx=\varepsilon^{n-\frac{2n}{n+1}}n^{\frac{2n}{n+1}}\int_{\mathrm{B}^+_q(\frac{\rho}{\varepsilon})}f^n dx\leq\int_{\mathbb{R}^n_+} f^ndx
\end{eqnarray*}

and:
 \begin{eqnarray*}
 \int_{\mathrm{B}^+_q(2\rho)\setminus\mathrm{B}^+_q(\rho)}|\D_g\psi^\pm_{\varepsilon}|^{\frac{2n}{n+1}}dx\leq C\varepsilon^{\frac{n(2n-1)}{n+1}}
 +C\varepsilon^{\frac{n(3n-1}{n+1}}\leq C\varepsilon^{\frac{n(2n-1)}{n+1}}.
\end{eqnarray*}

These estimates lead to:
\begin{eqnarray*}
\Big(\int_\M|\D_g\psi^\pm_{\varepsilon}|^{\frac{2n}{n+1}}dv(g)\Big)^{\frac{n+1}{n}}\leq\varepsilon^{n-1}n^2\mathrm{I}^{1+\frac{1}{n}}\Big(1+C\varepsilon^{\frac{n^2}{n+1}}\Big)=\varepsilon^{n-1}n^2\mathrm{I}^{1+\frac{1}{n}}\big(1+o(\varepsilon^{n-1})\big)
\end{eqnarray*}

where $\mathrm{I}=\int_{\mathbb{R}^n_+}f^ndx$. For the next, we let $\kappa=\mathrm{Re}\,\<\psi_0^\pm,\psi_1^\pm\>$ and we focus on an estimate of the denominator of the variational characterization~(\ref{charvar}) of $\lambda_{\min}(\M,\pa\M)$. We write:
\begin{eqnarray*}
\mathrm{Re}\<\D_g\psi^\pm_{\varepsilon},\psi^\pm_{\varepsilon}\>_{|\mathrm{B}^+_q(\rho)} & = & \mathrm{Re}\<\frac{n}{\varepsilon}f\big(\frac{r}{\varepsilon}\big)^{\frac{n}{2}+1}\big(1-\frac{x}{\varepsilon}\big)\cdot\psi_0^\pm,f\big(\frac{r}{\varepsilon}\big)^{\frac{n}{2}}\big(1-\frac{x}{\varepsilon}\big)\cdot\psi_0^\pm-\varepsilon_0\psi_1^\pm\>\\
 & = & \frac{n}{\varepsilon}f\big(\frac{r}{\varepsilon}\big)^n-\frac{n}{\varepsilon}\varepsilon_0\kappa\, f\big(\frac{r}{\varepsilon}\big)^{\frac{n}{2}+1}+\frac{n}{\varepsilon}\varepsilon_0 f\big(\frac{r}{\varepsilon}\big)^{\frac{n}{2}+1}\mathrm{Re}\<\big(\frac{x}{\varepsilon}\big)\cdot\psi_0^\pm,\psi_1^\pm\>.
\end{eqnarray*}

Integrating on $\mathrm{B}^+_q(\rho)$ gives:
\begin{eqnarray*}
\int_{\mathrm{B}^+_q(\rho)} \mathrm{Re}\<\D_g\psi^\pm_{\varepsilon},\psi^\pm_{\varepsilon}\> dx& = & \frac{n}{\varepsilon}\int_{\mathrm{B}^+_q(\rho)}f\big(\frac{r}{\varepsilon}\big)^n dx-\frac{n}{\varepsilon}\varepsilon_0\kappa\int_{\mathrm{B}^+_q(\rho)}f\big(\frac{r}{\varepsilon}\big)^{\frac{n}{2}+1}dx\\
& & +\frac{n}{\varepsilon}\varepsilon_0\int_{\mathrm{B}^+_q(\rho)}f\big(\frac{r}{\varepsilon}\big)^{\frac{n}{2}+1}\mathrm{Re}\<\big(\frac{x}{\varepsilon}\big)\cdot\psi_0^\pm,\psi_1^\pm\>dx\\
& = & n\varepsilon^{n-1}\Big(\int_{\mathrm{B}^+_q(\frac{\rho}{\varepsilon})}f(r)^n dx-\varepsilon_0\kappa\int_{\mathrm{B}^+_q(\frac{\rho}{\varepsilon})}f(r)^{\frac{n}{2}+1}dx+\mathrm{A}_{\varepsilon}\Big)
\end{eqnarray*}

where:
\begin{eqnarray*}
\mathrm{A}_\varepsilon=n\varepsilon_0\varepsilon^{-n}\int_{\mathrm{B}^+_q(\rho)}f\big(\frac{r}{\varepsilon}\big)^{\frac{n}{2}+1}\mathrm{Re}\<\frac{x}{\varepsilon}\cdot\psi_0^\pm,\psi_1^\pm\>dx.
\end{eqnarray*}

However $\mathrm{A}_\varepsilon=0$ since by hypothesis, we assumed that:
\begin{eqnarray*}
\mathrm{Re}\,\<\frac{x}{\varepsilon}\cdot\psi_0^\pm,\psi_1^\pm\>=0.
\end{eqnarray*}

Moreover an easy computation leads to:
\begin{eqnarray*}
\int_{\mathrm{B}^+_q(\frac{\rho}{\varepsilon})}f(r)^n dx=\mathrm{I}+O(\varepsilon^{\frac{n^2}{n+1}})
\end{eqnarray*}

and since $\varepsilon_0\sim\varepsilon^{n-1}$ when $\varepsilon\rightarrow 0$, we find:
\begin{eqnarray*}
\int_{\mathrm{B}^+_q(\rho)}\mathrm{Re}\<\D_g\psi^\pm_{\varepsilon},\psi^\pm_{\varepsilon}\> dx\geq n\varepsilon^{n-1}\Big(\mathrm{I}-\mathrm{C}_0\kappa\varepsilon^{n-1}+o(\varepsilon^{n-1})\Big)
\end{eqnarray*}

where $\mathrm{C}_0=\int_{\mathbb{R}^n_+}f(r)^{\frac{n}{2}+1}dx$. On the other hand, we compute:
\begin{eqnarray*}
\mathrm{Re}\<\D_g\psi_{\varepsilon}^\pm,\psi_{\varepsilon}^\pm\>_{|\mathrm{B}^+_q(2\rho)\setminus\mathrm{B}^+_q(\rho)} & = & -\mathrm{Re}\<\varepsilon_0\na\eta\cdot\theta^\pm,\varepsilon_0(\psi^\pm-\eta\theta^\pm)+\eta f\big(\frac{\rho}{\varepsilon}\big)^{\frac{n}{2}}\psi_0^\pm\>\\
& & +\mathrm{Re}\,\<f\big(\frac{\rho}{\varepsilon}\big)^{\frac{n}{2}}\na\eta\cdot\psi_0^\pm,\varepsilon_0(\psi^\pm-\eta\theta^\pm)+\eta f\big(\frac{\rho}{\varepsilon}\big)^{\frac{n}{2}}\psi_0^\pm\>\\
 & = & -\mathrm{Re}\<\varepsilon_0\na\eta\cdot\theta^\pm+f\big(\frac{\rho}{\varepsilon}\big)^{\frac{n}{2}}\na\eta\cdot\psi_0^\pm,\varepsilon_0\psi^\pm\>
\end{eqnarray*}

since $\mathrm{Re}\<\na\eta\cdot\theta^\pm,\theta^\pm\>=0$, $\mathrm{Re}\<\na\eta\cdot\psi_0^\pm,\psi_0^\pm\>=0$ and
\begin{eqnarray*}
 \mathrm{Re}\<\na\eta\cdot\psi_0^\pm,\theta^\pm\>=-\mathrm{Re}\<\na\eta\cdot\theta^\pm,\psi_0^\pm\>.
\end{eqnarray*}
 
This leads to the following estimation:
\begin{eqnarray*}
\mathrm{Re}\<\D_g\psi_{\varepsilon}^\pm,\psi_{\varepsilon}^\pm\>_{|\mathrm{B}^+_q(2\rho)\setminus\mathrm{B}^+_q(\rho)} \geq C\varepsilon^{2n-2}\rho^{1-n}
\end{eqnarray*}

and thus we have:
\begin{eqnarray*}
\int_{\mathrm{B}^+_q(2\rho)\setminus\mathrm{B}^+_q(\rho)}\mathrm{Re}\<\D_g\psi^\pm_{\varepsilon},\psi^\pm_{\varepsilon}\>dx=o(\varepsilon^{2(n-1)}).
\end{eqnarray*}

Now using the fact that $\mathrm{Re}\<\D_g\psi_{\varepsilon}^\pm,\psi_{\varepsilon}^\pm\>_{|\M\setminus\mathrm{B}^+_q(2\rho)}=0$, we get:
\begin{eqnarray*}
\int_{\M}\mathrm{Re}\<\D_g\psi_{\varepsilon}^\pm,\psi_{\varepsilon}^\pm\>\geq n\varepsilon^{n-1}\mathrm{I}\Big(1-\mathrm{C_0}\kappa\varepsilon^{n-1}++o(\varepsilon^{n-1})\Big).
\end{eqnarray*}

The variational characterization~(\ref{charvar}) of $\lambda_{\min}(\M,\pa\M)$ gives:
\begin{eqnarray*}
\lambda_{\min}(\M,\pa\M)\leq n\mathrm{I}^{\frac{1}{n}}\frac{1+o(\varepsilon^{n-1})}{1-\mathrm{C_0}\kappa\varepsilon^{n-1}+o(\varepsilon^{n-1})}=\lambda_{\min}(\hs,\pa\hs)+\mathrm{C}_0\kappa\varepsilon^{n-1}+o(\varepsilon^{n-1})
\end{eqnarray*}

where $\mathrm{C}_0$ is a positive constant. However, we assumed that $\kappa:=\mathrm{Re}\,\<\psi_0^\pm,\psi_1^\pm\><0$ so we can finally conclude that: 
\begin{eqnarray*}
\lambda_{\min}(\M,\pa\M)<\lambda_{\min}(\hs,\pa\hs).
\end{eqnarray*}

\hfill$\square$\\

The proof of Theorem~\ref{lcfm} is then reduced to prove the existence of a spinor field satisfying the hypothesis of Theorem~\ref{cop}.\\

{\it Proof of Theorem~\ref{lcfm}:}
We show that under the hypothesis of Theorem~\ref{lcfm}, there exists a spinor field on $\M\setminus\{q\}$ satisfying (\ref{hyp1}), (\ref{hyp2}) and (\ref{hyp3}). Without loss of generality, we can assume (using the conformal covariance of $\lambda_{\min}(\M,\pa\M)$) that for the metric $g$, there exists a boundary point $q\in\pa\M$ with locally flat and umbilic boundary neighbourhood. On the other hand, since the Dirac operator is supposed to be invertible on $\mathcal{H}^\pm$,  Proposition~\ref{Devfoncgreen} allows to write that (around a point $q\in\pa\M$), the chiral Green function admits the following development:
\begin{eqnarray*}
\G_\CHI^\pm(x,q)\psi_0^\pm=-\frac{2}{\omega_{n-1}}\frac{x-q}{|x-q|^n}\cdot\psi_0^\pm+\mathrm{m}_\CHI^\pm(x,q)\psi_0^\pm,
\end{eqnarray*}

where $\psi_0^\pm$ is a spinor such that $\nu\cdot\Ga\psi_0^\pm=\mp\psi_0^\pm$ and $\BCHI^\pm\big(\mathrm{m}_\CHI^\pm(\,.\,,q)\psi^\pm_{0|\pa\M}\big)=0$. Since the chiral mass operator is supposed to be non identically zero, we can choose a spinor $\psi_0^\pm$ which is an eigenspinor for the chiral mass operator $\mathrm{m}_\CHI^\pm(q)$ associated with the eigenvalue $\pm\kappa$ with $\kappa\in\mathbb{R}$. Moreover, Corollary~\ref{spectrum} insures that the eigenvalues $\pm\kappa$ are of opposite signs and so one of the two is positive and the other is negative. However the point $(1)$ of Remark~\ref{remark2} tells us that the chiral bag invariant $\lambda_{\min}(\M,\pa\M)$ does not depend of the choice of the boundary condition and we can then choose $\psi_0^-$ (for example) as being an eigenspinor for the negative chiral mass operator $\mathrm{m}^-_\CHI(q)$ associated with the eigenvalue $\kappa$ with $\kappa>0$. Then the chiral Green function $\G_\CHI^-(\,.\,,q)\psi_0^-$ is given by:
\begin{eqnarray*}
\G_\CHI^-(x,q)\psi_0^-=-\frac{2}{\omega_{n-1}}\frac{x-q}{|x-q|^n}\cdot\psi_0^-+\kappa\psi_0^-.
\end{eqnarray*}

Using the notations of Theorem~\ref{cop}, we have $\psi_1^-=-\kappa\psi_0^-$ and then we obtain:
\begin{eqnarray*}
\mathrm{Re}\,\<x\cdot\psi_0^-,\psi_1^-\>=\kappa\,\mathrm{Re}\,\<x\cdot\psi_0^-,\psi_0^-\>=0.
\end{eqnarray*}

Thus the spinor field $-\frac{\omega_{n-1}}{2}\G_\CHI^-(\,.\,,q)\psi_0^-$ satisfies the properties (\ref{hyp1}), (\ref{hyp2}) and (\ref{hyp3}) and so Theorem~\ref{cop} allows to conclude.
\hfill$\square$\\

As an application of this result, we obtain a spinorial proof of the Yambe problem on manifolds with boundary. In fact, we have:
\begin{corollary}
Suppose that $\mu(\M,\pa\M)>0$ and that the manifold $(\M^n,g,\sigma)$ ($n\geq 3)$ has a locally conformally flat with umbilic boundary point $q\in\pa\M$. Assume moreover that the chiral mass operator is non identically zero, then there exists a metric $\ov{g}\in [g]$ such that the scalar curvature $\mathrm{R}_{\ov{g}}$ is a positive constant and the mean curvature $\mathrm{H}_{\ov{g}}$ is zero.
\end{corollary}

{\it Proof:}
A sufficient condition for the existence of a solution of the Yamabe problem on manifolds with boundary is (see \cite{escobar:92}):
\begin{eqnarray}\label{solyam}
\mu(\M,\pa\M)<\mu(\hs,\pa\hs),
\end{eqnarray}

where $\mu(\M,\pa\M)$ is the Yamabe invariant of $\M$ endowed with the conformal class of $g$. The hypothesis $\mu(\M,\pa\M)>0$ implies that the Dirac operator is invertible, under the chiral bag boundary condition. Since the manifold has a locally conformally flat with umbilic boundary point $q\in\pa\M$ and since the chiral mass operator is supposed to be non zero, then we can apply Theorem~\ref{lcfm} and we finally get:
\begin{eqnarray*}
\lambda_{\min}(\M,\pa\M)<\lambda_{\min}(\hs,\pa\hs).
\end{eqnarray*}

However, using the Hijazi inequality~(\ref{hijbord}) and the fact that
\begin{eqnarray*}
 \lambda_{\min}(\hs,\pa\hs)=\frac{n}{2}\Big(\frac{\omega_{n}}{2}\Big)^{\frac{1}{n}}, 
\end{eqnarray*}
 
we obtain (\ref{solyam}) and thus the result follows direcly.
\hfill$\square$\\

This spinorial proof of the Yamabe problem problem on manifolds with boundary is obviously more difficult than the one given by Escobar. However, it seems that the inequality proved in Theorem~(\ref{lcfm}) gives existence of solutions to a Yamabe type problem for the Dirac operator under the chiral bag boundary condition involving critical Sobolev embedding theorems. A similar equation has been treated by B. Ammann for the boundaryless case (see \cite{habilbernd} or \cite{amm1}).


\section{The $\MIT$ Green function}


In this section, we construct the Green function for the Dirac operator under the $\MIT$ bag boundary condition. As an application of this construction, we give a simple proof of a positive mass theorem proved by Escobar \cite{escobar:92} in the context of the Yamabe problem on manifolds with boundary.


\subsection{Definitions and first properties}\label{definition1}


The $\MIT$ bag boundary condition has been introduced by physicists of the Massachusetts Institute of Technology in the seventies (see \cite{cjjtw}, \cite{cjjt} or \cite{johnson}). The spectrum of the Dirac operator under this boundary condition has then been studied on compact Riemannian spin manifolds with boundary (see \cite{hijazi.montiel.roldan:01} and \cite{sr1}). Another very nice application of this boundary condition can be found in \cite{hijazi.montiel.zhang:2} where the authors use its conformal covariance to give some estimates of the boundary Dirac operator involving a conformal invariant (which can be seen as an extrinsic version of the Hijazi inequality \cite{hijazi:86}). \\

We now briefly recall the definition of the $\MIT$ bag boundary condition. Consider the linear map
\begin{eqnarray*}
 i\nu\cdot:\mathbf{S}_g\longrightarrow\mathbf{S}_g
\end{eqnarray*}

which is an involution, so we can define the pointwise orthogonal projection
$$\begin{array}{lccl}
\BMIT^\pm : & \L^2(\mathbf{S}_g) & \longrightarrow & \L^2(\V^{\pm}_g)\\
 & \varphi & \longmapsto & \frac{1}{2}(\Id\pm i\nu\cdot)\varphi,
\end{array}$$ 

where  $\mathrm{V}_g^\pm$ is the eigensubbundle associated with the eigenvalues $\pm 1$ of the endomorphism $i\nu$. We can then check (see \cite{hijazi.montiel.roldan:01}) that these projections define elliptic boundary conditions for the Dirac operator. We can now define the Green function for the Dirac operator under the $\MIT$ bag boundary condition. 
\begin{definition}\label{MITgreen}
A Green function for the Dirac operator under the $\MIT$ bag boundary condition (or a $\MIT$ Green function) is given by a smooth section:
\begin{eqnarray}\label{MITgreen1}
\G^\pm_\MIT:\M\times\M\setminus\Delta\rightarrow\Sigma_g(\M)\boxtimes\big(\Sigma_g(\M)\big)^\ast
\end{eqnarray}

locally integrable on $\M\times\M$ which satisfies, in a weak sense, the following boundary problem:
$$\left\lbrace
\begin{array}{ll}
\D_g\big(\G^\pm_\CHI(x,y)\big)=\delta_y\mathrm{Id}_{\Sigma_y\M}\\ \\
\BMIT^\pm\big(\G^\pm_\CHI(x,y)\big)=0,\quad\textrm{for}\;x\in\pa\M\setminus\{y\},
\end{array}
\right.$$

for all $x\neq y\in\M$. In other words, we have:
\begin{eqnarray*}
\int_{\M}\<\G_{\MIT}^\pm(x,y)\psi_0^\pm,\D_g\varphi(x)\>dv(x)=\<\psi_0^\pm,\varphi(y)\>
\end{eqnarray*}

for all $y\in\M$, $\psi_0^\pm\in\Sigma_{y}(\M)$ satisfying $i\nu\cdot\psi_0^\pm=\pm\psi_0^\pm$ et $\varphi\in\Gamma\big(\Sigma_g(\M)\big)$ such that $\BMIT^\mp(\varphi_{|\pa\M})=0$.
\end{definition}
 
\begin{remark}
If we let:
\begin{eqnarray*}
\D_g^\pm:\mathcal{H}^\pm_g:=\{\phi\in\mathrm{H}^2_1(\Sigma\M)\;/\;\BMIT^\pm(\phi_{|\pa\M})=0\}\longrightarrow\L^2\big(\Sigma_g(\M)\big),
\end{eqnarray*}

we can easily check (using the Green formula) that $(\D^\pm_g)^\ast=\D^\mp_g$ where $(\D^\pm_g)^\ast$ is the formal adjoint of the Dirac operator under the boundary condition $\BMIT^\pm$. Thus the domain $\mathrm{dom}\,(\D^\pm_g)^\ast$ of the Dirac operator's adjoint is given by $\mathcal{H}^\mp_g$. The preceding definition of the $\MIT$ Green fonction is then consistent with the definition of weak solution of an equation. Indeed, the section $\G_\MIT^\pm(.\,,q)\psi_0^\pm\in\Gamma\big(\Sigma_g(\M)\big)$ satisfies on $\mathcal{H}^\pm_g$ and in a week sense the equation:
\begin{eqnarray*}
\D_g\big(\G^\pm_{\MIT}(.\,,q)\psi_0^\pm\big) = \delta_q
\end{eqnarray*}

if for all $\psi_0^\pm\in\Sigma_q\M$ such that $i\nu\cdot\psi_0^\pm=\pm\psi_0^\pm$:
\begin{eqnarray*}
\int_{\M}\<\G_{\MIT}(x,q)\psi_0^\pm,(\D^\pm_g)^\ast\phi\>dv(g)  = \<\psi_0^\pm,\phi(q)\>_{\Sigma_q\M}
\end{eqnarray*}

for all $\phi\in\mathrm{dom}\,(\D^\pm_g)^\ast=\mathcal{H}^\mp_g$. This definition is exactly the one given in Definition \ref{MITgreen}.
\end{remark}
 
In order to give a development of the $\MIT$ Green function in a neighbourhood of $q\in\pa\M$, it is very useful to study the behaviour of the $\MIT$ condition under the trivialization of the spinor bundle induced by Fermi coordinates and given in Section~\ref{GFDOCBC}. More precisely, we prove the following lemma:
\begin{lemma}\label{MITtri}
Let $\U$ and $\V$ be the open sets of the trivialization given in Section~\ref{GFDOCBC} and let $\Phi_0\in\Ga\big(\Sigma_{\xi}(\mathbb{R}^n_+)\big)$ be a parallel spinor such that
\begin{eqnarray*}
i\nu\cdot\ov{\Phi}_{0}(q)=\pm\ov{\Phi}_{0}(q)
\end{eqnarray*}

at one point $q\in\V\cap\pa\M$. Then we get:
\begin{eqnarray*}
i\nu\cdot\ov{\Phi}_{0|\V\cap\pa\M}=\pm\ov{\Phi}_{0|\V\cap\pa\M},
\end{eqnarray*}

that is $\BMIT^\mp(\ov{\Phi}_{0|\V\cap\pa\M})=0$.
\end{lemma}

\noindent {\it Proof:} One considers the function defined on $\V$ by $f(p)=|i\nu\cdot\ov{\Phi}_0-\ov{\Phi}_0|^2(p)$ and we then show that $f$ vanishes on $\V\cap\pa\M$. Indeed, for $1\leq i\leq n-1$, we have:
\begin{eqnarray*}
 e_i(f) & = & e_i(|i\nu\cdot\ov{\Phi}_0-\ov{\Phi}_0|^2) \\
 & = & 2\,e_i\big(|\ov{\Phi}_0|^2-i\<\nu\cdot\ov{\Phi}_0,\ov{\Phi}_0\>\big)\\
 & = & 2\,e_i\big(|\ov{\Phi}_0|^2+2\,\mathrm{Im}\<\nu\cdot\ov{\Phi}_0,\ov{\Phi}_0\>\big),
\end{eqnarray*}

where $\mathrm{Im}\,(z)$ is the imaginary part of a complex number $z\in\mathbb{C}$. However, since the spinor field $\Phi_0$ is parallel, we can assume that $|\Phi_0|^2=1$ and since the trivialization of the spinor bundle is an isometry, we get $|\ov{\Phi}_0|^2=1$ and so $e_i\big(|\ov{\Phi}_0|^2\big)=0$. Using the compatibilty of the spinorial Levi-Civita connection with the Hermitian metric leads to:
\begin{eqnarray*}
e_i\big(\mathrm{Im}\<\nu\cdot\ov{\Phi}_0,\ov{\Phi}_0\>\big)=\mathrm{Im}\<\ov{\na}_{e_i}\nu\cdot\ov{\Phi}_0,\ov{\Phi}_0\>+2\,\mathrm{Im}\<\nu\cdot\ov{\na}_{e_i}\ov{\Phi}_0,\ov{\Phi}_0\>.
\end{eqnarray*}

The local expression of the spinorial Levi-Civita connection $\ov{\na}$ acting on $\Sigma_g(\M)$ gives:
\begin{eqnarray*}
\mathrm{Im}\<\nu\cdot\ov{\na}_{e_i}\ov{\Phi}_0,\ov{\Phi}_0\> & = & \frac{1}{4}\sum_{1\leq j\neq k\leq n-1}\widetilde{\Ga}_{ij}^k\,\mathrm{Im}\<\nu\cdot e_j\cdot e_k\cdot\ov{\Phi}_0,\ov{\Phi}_0\>\\ 
& & -\frac{1}{2}\,\mathrm{Im}\<\ov{\na}_{e_i}\nu\cdot\ov{\Phi}_0,\ov{\Phi}_0\>.
\end{eqnarray*} 

Note that:
\begin{eqnarray*}
\<\nu\cdot e_j\cdot e_k\cdot\ov{\Phi}_0,\ov{\Phi}_0\> & = & \<\ov{\Phi}_0,\nu\cdot e_j\cdot e_k\cdot\ov{\Phi}_0\>,
\end{eqnarray*}

and so we have $e_i(f)=0$ for all $1\leq i\leq n-1$. Since we assumed that $f(q)=0$, we obtain immediately the result of this lemma.
\hfill$\square$\\

We can now state the analogous result of Proposition~\ref{Devfoncgreen} for the $\MIT$ Green function. Indeed, we have:
\begin{proposition}\label{DevfoncgreenMIT}
Assume that there exists a point $q\in\pa\M$ which has a locally flat with umbilic boundary neighbourhood. Then the $\MIT$ Green function exists and admits the following development near $q$:
\begin{eqnarray}\label{devfoncgreenmit}
\G_\MIT^\pm(x,q)\psi_0^\pm=-\frac{2}{\omega_{n-1}}\frac{x-q}{|x-q|^n}\cdot\psi_0^\pm+\mathrm{m}_\MIT^\pm(x,q)\psi_0^\pm,
\end{eqnarray}

where $\psi_0^\pm\in\Sigma_q(\M)$ is a spinor field satisfying $i\,\nu\cdot\psi_0^\pm=\pm\psi_0^\pm$ (in other words $\BMIT^{\mp}(\psi_0^\pm)=0$) and where $\mathrm{m}_\MIT^\pm(\,.\,,q)\psi_0^\pm$ is a smooth spinor field such that $\D_g\big(\mathrm{m}_\MIT^\pm(\,.\,,q)\psi_0^\pm\big)=0$ around $q$. Moreover, along $\pa\M$, we have:
\begin{eqnarray}\label{massmit}
\BMIT^{\mp}\big(\mathrm{m}_\MIT^\pm(.\,,q)\psi^\pm_{0|\pa\M}\big)=0.
\end{eqnarray}
\end{proposition}
 
{\it Proof:} The proof of this proposition follows the proof of Proposition~\ref{Devfoncgreen} using Lemma~\ref{MITtri} and the fact that the Dirac operator:
\begin{eqnarray*}
\D_g:\mathcal{H}^\pm_g:=\{\phi\in\mathrm{H}^2_1(\Sigma\M)\;/\;\BMIT^\pm(\phi_{|\pa\M})=0\}\longrightarrow\L^2\big(\Sigma_g(\M)\big),
\end{eqnarray*}

defines an invertible operator. 
\hfill$\square$\\


\subsection{Application: The positive mass theorem for the Yamabe problem on manifolds with boundary}


First we recall briefly the notion of mass for compact manifolds with boundary. For more details, we refer to \cite{lee.parker:87} for the boundaryless case and \cite{escobar:92} for the non empty boundary case. Let $(\M^n,g)$ be a connected compact Riemannian manifold which is locally conformally flat with umbilic boundary. If we assume that the Yamabe invariant $\mu(\M,\pa\M)$ of the manifold $\M$ is positive and if $q\in\pa\M$ then there exists a unique Green function $\mathcal{G}_q$ for the conformal Laplacian 
\begin{eqnarray*}
\L_g:=4\frac{n-1}{n-2}\Delta_g+\R_g
\end{eqnarray*}
 
under the boundary condition 
\begin{eqnarray*}
\mathrm{B}_g:= -\frac{2}{n-2}\frac{\pa}{\pa\nu}+\H_g,
\end{eqnarray*}

that is, a smooth function $\mathcal{G}_q$ defined on $\M\setminus\{q\}$ which satisfies, in a weak sense, the boundary problem:
\begin{equation}
\left\lbrace
\begin{array}{lll}\label{GFCF}
\L_g\mathcal{G}_q & = & \delta_q \quad\text{on}\;\M \\ \\
\mathrm{B}_g\mathcal{G}_{q|\pa\M} & = & \delta_q \quad\text{along}\;\pa\M
\end{array}
\right.
\end{equation}

One can then check that if there exists a point $q\in\pa\M$ with a locally flat and umbilic boundary neighbourhood in the metric $\ov{g}\in [g]$, then the Green function $\mathcal{G}_q$ admits the following expansion (near $q$):
\begin{eqnarray}\label{devfongrlc}
\mathcal{G}_q(x)=\frac{1}{(n-2)\omega_{n-1} r^{n-2}}+\mathrm{A}+\alpha_{q}(x)
\end{eqnarray}

where $\mathrm{A}\in\mathbb{R}$ and $\alpha_q$ is a harmonic function (near $q$) which satisfies $\alpha_q(q)=0$ and $\frac{\pa\alpha_q}{\pa\nu}=0$. The positive mass theorem proved by Escobar can then be stated as follows:
\begin{center}
{\it The constant $\mathrm{A}$ satisfies $\mathrm{A}\geq 0$. Moreover, $\mathrm{A}=0$ if and only if $\M$ is conformally isometric to the round hemisphere $\mathbb{S}^n_+$.\\}
\end{center}

We give here a proof of this result for spin manifolds. However, in our proof, it is not necessary to impose that the manifold is locally conformally flat with umbilic boundary. Indeed, it is sufficient to assume that the manifold has only one point $q\in\pa\M$ which has a locally conformally flat with umbilic boundary neighbourhood. 

The proof of this positive mass theorem is inspired of the work of Ammann and Humbert \cite{amm2} and is based on the construction of the Green function of the Dirac operator. The chiral Green function used in Section~\ref{GFDOCBC} seems to be a good candidate for our purpose however it needs the existence of a chirality operator which is not the case in any dimension. That is why we will use another conformal covariant boundary condition which exists without additional assumptions, the $\MIT$ bag boundary condition. More precisely, we prove:
\begin{theorem}
Let $(\M^n,g,\sigma)$ be a connected compact Riemannian spin manifold with non empty smooth boundary. Assume that the Yamabe invariant is positive and that there exists a point $q\in\pa\M$ which has a locally conformally flat with umbilic boundary neighbourhood. Then the mass ${\rm A}$ of the manifold $\M$ satisfies ${\rm A}\geq 0$. Moreover, equality holds if and only if $\M$ is conformally isometric to the round hemisphere $\hs$.
\end{theorem}

{\it Proof:}
First note that we can assume that for the metric $g$, the point $q\in\pa\M$ has a locally flat with umbilic boundary neighbourhood. Proposition~\ref{DevfoncgreenMIT} allows to deduce that there exists a unique $\MIT$ Green function $\G_{\MIT}^-(\,.\,,q)$ for the Dirac operator which admits the following development near $q\in\pa\M$:
\begin{eqnarray}
\G_{\MIT}^-(x,q)\psi_0 & = & \G_{\eucl}(x,q)\psi_0+\mathrm{m}^-_{\MIT}(x,q)\psi_0,
\end{eqnarray}

where $\psi_0\in\Sigma_q\M$ is such that $i\nu\cdot\psi_0=-\psi_0$ and $\mathrm{m}_\MIT^-(\,.\,,q)\psi_0\in\Gamma\big(\Sigma_g(\M)\big)$ is a harmonic spinor near $q$ which satisfies:
\begin{eqnarray*}
\BMIT^-\big(\mathrm{m}_\MIT^-(.\,,q)\psi_{0|\pa\M}\big)=0.
\end{eqnarray*}

Without loss of generality, we can assume that $|\psi_0|^2=1$. On the other hand, since $\mu(\M,\pa\M)>0$, there exists a unique Green function $\mathcal{G}_q$ for the conformal Laplacian smooth on $\M\setminus\{q\}$ such that $\mathcal{G}_q>0$ (see \cite{escobar:92}). Moreover, since the point $q$ is supposed to have a locally flat with umbilic boundary neighbourhood in the metric $g$, the Green function $\mathcal{G}_q$ can be written 
\begin{eqnarray*}
\mathcal{G}_q(x)=\frac{1}{(n-2)\omega_{n-1} r^{n-2}}+\mathrm{A}+\alpha_{q}(x)
\end{eqnarray*}

near $q\in\pa\M$ and where $\mathrm{A}\in\mathbb{R}$ and $\alpha_q$ is a smooth function on $\M$. Now consider the conformal change of the metric on $\M\setminus\{q\}$ given by:
 \begin{eqnarray*}
 \ov{g}=\big((n-2)\omega_{n-1}\mathcal{G}_q\big)^{\frac{4}{n-2}}g =\widetilde{\mathcal{G}}_q^{\frac{4}{n-2}}g.
\end{eqnarray*}

Since the Green function satisfies the boundary problem~(\ref{GFCF}) and using the fact that the scalar and mean curvatures in $g$ and $\ov{g}$ are related by:
\begin{eqnarray*}
\ov{\R}=\widetilde{\mathcal{G}}_q^{-\frac{n+2}{n-2}}\L_g\widetilde{\mathcal{G}}_q\quad\textrm{and}\quad\ov{\H}=\widetilde{\mathcal{G}}_q^{-\frac{n}{n-2}}\mathrm{B}_g\widetilde{\mathcal{G}}_q,
\end{eqnarray*}

we obtain that the scalar curvature of $(\M\setminus\{q\},\ov{g})$ is $\ov{\R}=0$ and the mean curvature of $(\pa\M\setminus\{q\},\ov{g})$ in $\M$ is $\ov{\H}=0$. We can then identify the spinor bundle $(\M\setminus\{q\},g)$ with the one over $(\M\setminus\{q\},\ov{g})$ thanks to the bundle ismorphism $\F$ used in the proof of Proposition~\ref{Confgreen}. Now consider the spinor field defined by:
\begin{eqnarray*}
\ov{\G}_{q}(\,.\,)\ov{\psi}_0 = \widetilde{\mathcal{G}}_q^{-\frac{n+1}{n-2}}\F\big(\G_\MIT^-(\,.\,,q)\psi_0\big)\in\Gamma\big(\Sigma_{\ov{g}}(\M\setminus\{q\})\big).
\end{eqnarray*}

Using the conformal covariance (\ref{covconf}) of the Dirac operator, we have:
\begin{eqnarray}\label{grendi}
\D_{\ov{g}}\big(\ov{\G}_{q}(\,.\,)\ov{\psi}_0\big)=\widetilde{\mathcal{G}}_q^{-\frac{n-1}{n-2}}\F\Big(\D_g\big(\G_\MIT^-(\,.\,,q)\psi_0\big)\Big)=0
\end{eqnarray}

on $\M\setminus\{q\}$. Moreover since the $\MIT$ condition is also conformally invariant, we obtain:
\begin{eqnarray*}
\ov{\mathbb{B}}_\MIT^-\big(\ov{\G}_{q}(\,.\,)\ov{\psi}_{0|\pa\M}\big)=0,
\end{eqnarray*}

where $\ov{\mathbb{B}}_\MIT^-=\frac{1}{2}\big(\Id-i\ov{\nu}\,\ov{\cdot}\big)$ is the projection of the $\MIT$ bag boundary condition in the metric $\ov{g}$. Now let $\varepsilon>0$ and denote by $\mathrm{B}^+_q(\varepsilon)$ the half-ball centered at $q$ with radius $\varepsilon$. Using the formula~(\ref{grendi}), the Schr\"odinger-Lichnerowicz formula and the fact that the scalar curvature of $\M\setminus\{q\}$ is zero, we get:
\begin{eqnarray*}
0=\int_{\M\setminus\mathrm{B}^+_q(\varepsilon)}\<\D^2_{\ov{g}}\big(\ov{\G}_{q}(x)\ov{\psi}_0\big),\ov{\G}_{q}(x)\ov{\psi}_0\>dv(\ov{g})=\int_{\M\setminus\mathrm{B}^+_q(\varepsilon)}\<\ov{\na}^\ast\ov{\na}(\ov{\G}_{q}(x)\ov{\psi}_0),\ov{\G}_{q}(x)\ov{\psi}_0\>dv(\ov{g}).
\end{eqnarray*}

An integration by parts leads to:
\begin{eqnarray}\label{mass}
\int_{\pa\mathrm{B}^+_q(\varepsilon)}\<\ov{\na}_{\ov{\nu}_{\varepsilon}}(\ov{\G}_{q}(x)\ov{\psi}_0),\ov{\G}_{q}(x)\ov{\psi}_0\>ds(\ov{g})& = & \int_{\pa\M_\varepsilon}\<\ov{\na}_{\ov{\nu}}(\ov{\G}_{q}(x)\ov{\psi}_0),\ov{\G}_{q}(x)\ov{\psi}_0\>ds(\ov{g})\nonumber\\
& & + \int_{\M\setminus\mathrm{B}^+_q(\varepsilon)}|\ov{\na}(\ov{\G}_{q}(x)\ov{\psi}_0)|^2 dv(\ov{g})
\end{eqnarray}

where $\pa\M_{\varepsilon}=\pa\M\setminus\big(\pa\M\cap\pa\mathrm{B}^+_q(\varepsilon)\big)$ and $\nu$ (resp. $\nu_{\varepsilon}$) is the inner unit vector field (resp. outer) normal to $\pa\M_\varepsilon$ (resp. $\pa\mathrm{B}^+_q(\varepsilon)$) in the metric $g$. An easy calculation gives:
\begin{eqnarray*}
\ov{\na}_{\ov{\nu}}(\ov{\G}_{q}(x)\ov{\psi}_0) & = & \frac{n-1}{2}\ov{\H}\,\ov{\G}_{q}(x)\ov{\psi}_0-\ov{\nu}\,\ov{\cdot}\D_{\ov{g}}\big(\ov{\G}_{q}(x)\ov{\psi}_0\big)-\D^{\pa\M}_{\ov{g}}\big(\ov{\G}_{q}(x)\ov{\psi}_0\big)
\end{eqnarray*}

where 
\begin{eqnarray*}
\D^{\pa\M}_{g}:=\sum_{i=1}^{n-1} e_i\cdot\nu\cdot\na^{\bf{S}}_{e_i}
\end{eqnarray*}

is the boundary Dirac operator acting on the restricted spinor bundle over $\pa\M$ endowed with the metric $g$. However, on $\pa\M_\varepsilon$ we have $\ov{\H}=0$ and $\D_{\ov{g}}\big(\ov{\G}_{q}(x)\ov{\psi}_0\big)=0$, so (\ref{mass}) gives:
\begin{eqnarray*}
\int_{\pa\mathrm{B}^+_q(\varepsilon)}\<\ov{\na}_{\ov{\nu}_{\varepsilon}}(\ov{\G}_{q}(x)\ov{\psi}_0),\ov{\G}_{q}(x)\ov{\psi}_0\>ds(\ov{g})& = & -\int_{\pa\M_\varepsilon}\<\D^{\pa\M}_{\ov{g}}\big(\ov{\G}_{q}(x)\ov{\psi}_0\big),\ov{\G}_{q}(x)\ov{\psi}_0\>ds(\ov{g})\\
& & + \int_{\M\setminus\mathrm{B}^+_q(\varepsilon)}|\ov{\na}(\ov{\G}_{q}(x)\ov{\psi}_0)|^2 dv(\ov{g}).
\end{eqnarray*}

On the other hand, using the conformal covariance of the $\MIT$ bag boundary condition, we have:
\begin{eqnarray*}
i\ov{\nu}\,\ov{\cdot}\ov{\G}_{q}(x)\ov{\psi}_0=\ov{\G}_{q}(x)\ov{\psi}_0,
\end{eqnarray*}

for all $x\in\pa\M_\varepsilon$ and so:
\begin{eqnarray*}
\<\D^{\pa\M}_{\ov{g}}\big(\ov{\G}_{q}(x)\ov{\psi}_0\big),\ov{\G}_{q}(x)\ov{\psi}_0\> & = & \<i\ov{\nu}\,\ov{\cdot}\D^{\pa\M}_{\ov{g}}\big(\ov{\G}_{q}(x)\ov{\psi}_0\big),
i\ov{\nu}\,\ov{\cdot}\ov{\G}_{q}(x)\ov{\psi}_0\>\\
& = & -\<\D^{\pa\M}_{\ov{g}}\big(\ov{\G}_{q}(x)\ov{\psi}_0\big),\ov{\G}_{q}(x)\ov{\psi}_0\>.
\end{eqnarray*}

Hence Inequality~(\ref{mass}) gives:
\begin{eqnarray}
0 \leq \int_{\M\setminus\mathrm{B}^+_q(\varepsilon)}|\ov{\na}(\ov{\G}_{q}(x)\ov{\psi}_0)|^2 dv(\ov{g}) =  \frac{1}{2}\int_{\pa\mathrm{B}^+_q(\varepsilon)}\ov{\nu}_{\varepsilon}|\ov{\G}_{q}(x)\ov{\psi}_0|^2ds(\ov{g}).
\end{eqnarray}

The vector field $\ov{\nu}_\varepsilon$ is the inner normal field to $\pa\mathrm{B}^+_q(\varepsilon)$ for the metric $\ov{g}$ and $\ov{g}$ is conformal to $g$ (which is flat around $q$), then the vector field $\ov{\nu}_{\varepsilon}$ is colinear to $\frac{\pa}{\pa r}$, that is there exists a constant $c>0$ such that $\ov{\nu}_{\varepsilon}=-c\frac{\pa}{\pa r}$. However:
\begin{eqnarray*}
1=\ov{g}(\ov{\nu}_{\varepsilon},\ov{\nu}_{\varepsilon})=c^2\widetilde{\mathcal{G}}_q^{\frac{4}{n-2}}g(\frac{\pa}{\pa r},\frac{\pa}{\pa r})=c^2(n-2)^{\frac{4}{n-2}}\,\omega_{n-1}^{\frac{4}{n-1}}\,\mathcal{G}_q^{\frac{4}{n-2}} 
\end{eqnarray*}

and since the half-ball $\mathrm{B}_q^+(\varepsilon)$ is contained in the open (flat) set of the trivialization near $q$, the Green function for the conformal Laplacian admits the expansion (\ref{devfongrlc}), that is:
\begin{eqnarray*}
c^2(n-2)^{\frac{4}{n-2}}\,\omega_{n-1}^{\frac{4}{n-1}}\,\Big(\frac{1}{(n-2)\omega_{n-1} r^{n-2}}+\mathrm{A}+\alpha_{q}(x)\Big)^{\frac{4}{n-2}}=1.
\end{eqnarray*}

An easy calculation then shows that $c=\varepsilon^2+o(\varepsilon^2)$ and finally we have $\ov{\nu}_{\varepsilon}=-(\varepsilon^2+o(\varepsilon^2))\frac{\pa}{\pa r}$. We now give an estimate of  $\ov{\nu}_{\varepsilon}|\ov{\G}_{q}(x)\ov{\psi}_0|^2$ on $\pa\mathrm{B}_q^+(\varepsilon)$; for this, we write:
\begin{eqnarray*}
|\ov{\G}_{q}(x)\ov{\psi}_0|^2 & = & \widetilde{\mathcal{G}}_q^{-2\frac{n-1}{n-2}}|\G_\MIT^-(x,q)\psi_0|^2\\
& = & (n-2)^{-2\frac{n-1}{n-2}}\,\omega_{n-1}^{-2\frac{n-1}{n-2}}\mathcal{G}_q^{-2\frac{n-1}{n-2}}\,|\G_\MIT^-(x,q)\psi_0|^2\\ 
& = & \Big(\frac{1}{r^{n-2}}+(n-2)\omega_{n-1}\mathrm{A}+\widetilde{\alpha}_{q}(x)\Big)^{-2\frac{n-1}{n-2}}\,|\G_{\eucl}(x,q)\psi_0+\mathrm{m}^-_{\MIT}(x,q)\psi_0|^2
\end{eqnarray*}

where $\widetilde{\alpha}_q(x)=(n-2)\,\omega_{n-1}\alpha_{q}(x)$. Using the expansion~(\ref{MITgreen1}) of the $\MIT$ Green function, we get:
\begin{eqnarray*}
|\ov{\G}_{q}(x)\ov{\psi}_0|^2 & = & \Big(1+(n-2)\omega_{n-1}\mathrm{A}r^{n-2}+\widetilde{\alpha}_{q}(x)r^{n-2}\Big)^{-2\frac{n-1}{n-2}}\\
& & \times\Big(1+|\mathrm{m}^-_{\MIT}(x,q)\psi_0|^2+2\mathrm{Re}\,\<\G_{\eucl}(x,q)\psi_0,\mathrm{m}^-_{\MIT}(x,q)\psi_0\>\Big).
\end{eqnarray*}

Now note that, with some calculations, we obtain:
\begin{eqnarray*}
\frac{\pa}{\pa r}|\ov{\G}_{q}(x)\ov{\psi}_0|^2=-2(n-1)(n-2)\,\omega_{n-1}\,\mathrm{A}\,\varepsilon^{n-3}+o(\varepsilon^{n-3}).
\end{eqnarray*}

On the other hand, we have:
\begin{eqnarray*}
ds (\ov{g}) & = & \sqrt{\mathrm{det}(\ov{g})}dx\\
 & = & \widetilde{\mathcal{G}}_q^{2\frac{n-1}{n-2}}dx\\
 & = & \varepsilon^{-2(n-1)}\big(1+o(1)\big)dx
\end{eqnarray*}

where $dx$ is the standard volume form of $\mathbb{R}^{n-1}$. We then have finally shown that:
\begin{eqnarray*}
0\leq \int_{\M\setminus\mathrm{B}^+_q(\varepsilon)}|\ov{\na}(\ov{\G}_{q}(x)\ov{\psi}_0)|^2 dv(\ov{g}) & = & \frac{1}{2}\int_{\pa\mathrm{B}^+_q(\varepsilon)}\ov{\nu}_{\varepsilon}|\ov{\G}_{q}(x)\ov{\psi}_0|^2ds(\ov{g})\\
& = & \frac{1}{2}\mathrm{vol}\big(\pa\mathrm{B}^+_q(\varepsilon)\big)\big(-\varepsilon^2+o(\varepsilon^2)\big)\big(\varepsilon^{-2(n-1)}+o(\varepsilon^{-2(n-1)})\big)\\
& & \times \big(-2(n-1)(n-2)\omega_{n-1}\,\mathrm{A}\varepsilon^{n-3}+o(\varepsilon^{n-3})\big)\\
& = & (n-1)(n-2)\,\omega_{n-1}^2\,\mathrm{A}+o(1),
\end{eqnarray*}

hence $\mathrm{A}\geq0$. Now assume that $\mathrm{A}=0$; using the preceding identity, we have that the spinor field $\ov{\G}_{q}(\,.\,)\ov{\psi}_0$ is parallel on $\M\setminus\{q\}$, that is $\ov{\na}_X(\ov{\G}_{q}(\,.\,)\ov{\psi}_0)=0$ for all $X\in\Ga\big(\mathrm{T}(\M\setminus\{q\})\big)$. Since the choice of $\psi_0$ is arbitrary, we easily construct a basis of parallel spinor fields over $(\M\setminus\{q\},\ov{g})$. On the other hand, the spinor $\ov{\G}_{q}(\,.\,)\ov{\psi}_0$ satisfies the $\MIT$ bag boundary condition, i.e.:
\begin{eqnarray*}
i\ov{\nu}\,\ov{\cdot}\ov{\G}_{q}(\,.\,)\ov{\psi}_{0|\pa\M}=\ov{\G}_{q}(\,.\,)\ov{\psi}_{0|\pa\M}.
\end{eqnarray*}

So if we derive (along $\pa\M$) this equality, we get:
\begin{eqnarray*}
\ov{\na}_X\big(i\ov{\nu}\,\ov{\cdot}\ov{\G}_{q}(\,.\,)\ov{\psi}_0\big) & = & -i\ov{\mathrm{A}}(X)\ov{\cdot}\ov{\G}_{q}(\,.\,)\ov{\psi}_0=0
\end{eqnarray*}

since the spinor field $\ov{\G}_{q}(\,.\,)\ov{\psi}_0$ is parallel. Using the fact the spinor field $\ov{\G}_{q}(\,.\,)\ov{\psi}_0$ has no zero (since it is a constant section), we deduce that the boundary $\pa\M$ is totally geodesic in $\M$. The restriction to the boundary of a parallel spinor field from a basis of the spinor bundle over $(\M,\ov{g})$ gives a parallel spinor field on the restricted spinor bundle along the boundary. This is an easy consequence of the spinorial Gauss formula and of the fact that the boundary is totally geodesic. We can thus easily conclude that the boundary is isometric to the Euclidean space $\mathbb{R}^{n-1}$, hence the manifold $(\M\setminus\{q\},\ov{g})$ is isometric to the Euclidean half-space $(\mathbb{R}^n_+,\xi)$. Now consider $\mathrm{I}:(\M\setminus\{q\},\ov{g})\rightarrow(\mathbb{R}^n_+,\xi)$ an isometry and let $f(x)=1+\xi\big(f(x)\big)/4$ where $\xi\big(f(x)\big):=\xi\big(f(x),f(x)\big)$. Then $\M\setminus\{q\}$ endowed with the metric $f^{-2}\ov{g}=f^{-2}\widetilde{\mathcal{G}}_q^{\frac{4}{n-2}}g$ is isometric to the standard hemisphere $\hs$. The function $f^{-2}\widetilde{\mathcal{G}}_q^{\frac{4}{n-2}}$ is smooth on $\M\setminus\{q\}$ and extends to a positive function all over $\M$. Thus we have shown that $\M$ is conformally equivalent to the standard hemisphere $(\hs,g_{\rm st})$.
\hfill$\square$


\bibliographystyle{amsalpha}     
\bibliography{bibthese1}


\end{document}